\ifdefined\IcmlBuild
  \documentclass{article}
\else
  \documentclass[11pt]{article}
\fi
\newif\ificmlbuild
\ifdefined\IcmlBuild\icmlbuildtrue\else\icmlbuildfalse\fi
\newcommand{\IfIcmlBuild}[1]{\ificmlbuild#1\fi}
\newcommand{\IfNotIcmlBuild}[1]{\ificmlbuild\else#1\fi}

\ificmlbuild
  \usepackage{hyperref}
  \usepackage{icml2027}\else
  \usepackage[letterpaper,margin=1in]{geometry}
  \usepackage[square,comma,sort&compress]{natbib}
\fi
\usepackage{amsmath,amsthm,mathtools}
\usepackage{amssymb}
\ificmlbuild\else
  \usepackage{lmodern}
  \usepackage[T1]{fontenc}
  \usepackage{charter}
\fi
\usepackage{graphicx,booktabs,tabularx,xcolor}
\ificmlbuild\else
  \usepackage[bookmarksnumbered,bookmarksopen,unicode]{hyperref}
\fi
\usepackage[algo2e,ruled,vlined]{algorithm2e}

\SetCommentSty{algcommentfont}
\SetKwComment{Comment}{$\triangleright$\ }{}
\IfNotIcmlBuild{\widowpenalty=10000 \clubpenalty=10000 }\IfNotIcmlBuild{}
\newcolumntype{Y}{>{\raggedright\arraybackslash}X}
\definecolor{paperdarkblue}{HTML}{17365D}
\definecolor{leangreen}{HTML}{2E7D32}
\definecolor{verifgray}{HTML}{4A4A4A}
\hypersetup{colorlinks=true,allcolors=paperdarkblue}
\providecommand{\doi}{}
\renewcommand{\doi}[1]{doi: \href{https://doi.org/\detokenize{#1}}{\begingroup\urlstyle{rm}\nolinkurl{#1}\endgroup}}
\graphicspath{{images/}}

\newcommand{\R}{\mathbb R}
\newcommand{\E}{\mathbb E}

\newcommand{\eps}{\varepsilon}
\newcommand{\diam}{\operatorname{diam}}
\newcommand{\dist}{\operatorname{dist}}
\newcommand{\conv}{\operatorname{conv}}

\newcommand{\Id}{\mathrm{Id}}
\newcommand{\LMO}{\mathtt{LMO}}
\newcommand{\err}{\operatorname{err}}
\newcommand{\one}{\mathbf 1}
\newcommand{\cls}{\mathcal C}
\newcommand{\ip}[2]{\langle #1,#2\rangle}
\newcommand{\norm}[1]{\lVert #1\rVert}
\DeclareMathOperator*{\argmin}{arg\,min}

\theoremstyle{plain}
\newtheorem{theorem}{Theorem}[section]
\newtheorem{lemma}[theorem]{Lemma}
\newtheorem{proposition}[theorem]{Proposition}
\newtheorem{corollary}[theorem]{Corollary}
\theoremstyle{definition}
\newtheorem{definition}[theorem]{Definition}
\theoremstyle{remark}
\newtheorem{remark}[theorem]{Remark}
\newenvironment{proofsketch}[1][Proof sketch]{\begin{proof}[#1]}{\end{proof}}

\newif\ifusedeferredproofs
\ifdefined\EnableDeferredProofs\usedeferredproofstrue\else\usedeferredproofsfalse\fi
\ifusedeferredproofs
  \usepackage[createShortEnv,disablePatchSection]{proof-at-the-end}
  \NewDocumentEnvironment{deferredtheorem}{o +b}{\IfNoValueTF{#1}
      {\begin{theoremE}[][proof at the end, restate, category=deferred]}{\begin{theoremE}[#1][proof at the end, restate, category=deferred]}#2\end{theoremE}}{}
  \NewDocumentEnvironment{deferredlemma}{o +b}{\IfNoValueTF{#1}
      {\begin{lemmaE}[][proof at the end, restate, category=deferred]}{\begin{lemmaE}[#1][proof at the end, restate, category=deferred]}#2\end{lemmaE}}{}
  \NewDocumentEnvironment{deferredproposition}{o +b}{\IfNoValueTF{#1}
      {\begin{propositionE}[][proof at the end, restate, category=deferred]}{\begin{propositionE}[#1][proof at the end, restate, category=deferred]}#2\end{propositionE}}{}
  \NewDocumentEnvironment{deferredcorollary}{o +b}{\IfNoValueTF{#1}
      {\begin{corollaryE}[][proof at the end, restate, category=deferred]}{\begin{corollaryE}[#1][proof at the end, restate, category=deferred]}#2\end{corollaryE}}{}
  \NewDocumentEnvironment{deferredproof}{+b}{\begin{textAtEnd}[category=deferred]
\begin{proof}
#1
\end{proof}
      \end{textAtEnd}
  }{}
  \newcommand{\PrintDeferredProofs}{\IfFileExists{\pratendGeneratePrefixFile{\jobname}deferred.tex}{\section{Deferred proofs}\label{app:deferred}\printProofs[deferred]}{}}
\else
  \NewDocumentEnvironment{deferredtheorem}{o +b}{\IfNoValueTF{#1}{\begin{theorem}}{\begin{theorem}[#1]}#2\end{theorem}}{}
  \NewDocumentEnvironment{deferredlemma}{o +b}{\IfNoValueTF{#1}{\begin{lemma}}{\begin{lemma}[#1]}#2\end{lemma}}{}
  \NewDocumentEnvironment{deferredproposition}{o +b}{\IfNoValueTF{#1}{\begin{proposition}}{\begin{proposition}[#1]}#2\end{proposition}}{}
  \NewDocumentEnvironment{deferredcorollary}{o +b}{\IfNoValueTF{#1}{\begin{corollary}}{\begin{corollary}[#1]}#2\end{corollary}}{}
  \NewDocumentEnvironment{deferredproof}{+b}{\begin{proof}#1\end{proof}}{}
  \newcommand{\PrintDeferredProofs}{}
\fi

\newif\ifshowchanges
\showchangesfalse

\NewDocumentEnvironment{changedblock}{+b}{\ifshowchanges\begingroup\color{blue}#1\endgroup\else#1\fi}{}

\ificmlbuild
  \newenvironment{widefigure}[1][t]{\begin{figure*}[#1]}{\end{figure*}}
  \newenvironment{widetable}[1][t]{\begin{table*}[#1]}{\end{table*}}
\else
  \newenvironment{widefigure}[1][t]{\begin{figure}[#1]}{\end{figure}}
  \newenvironment{widetable}[1][t]{\begin{table}[#1]}{\end{table}}
\fi

\newcommand{\leanverified}{\texorpdfstring{\unskip\nobreak\hspace{0.4em}\textcolor{leangreen}{\textnormal{\textsf{\footnotesize Lean}\,$\checkmark$}}}{}}

\ificmlbuild
\icmltitlerunning{Smoothness of the Feasible Set Does Not Accelerate Linear-Oracle Methods}
\else
\hypersetup{pdfauthor={Sebastian Pokutta},
 pdftitle={Smoothness of the Feasible Set Does Not Accelerate Linear-Oracle Methods}}
\title{Smoothness of the Feasible Set Does Not Accelerate\\ Linear-Oracle Methods}
\author{Sebastian Pokutta\thanks{ORCID \href{https://orcid.org/0000-0001-7365-3000}{0000-0001-7365-3000}.}\\
 Zuse Institute Berlin and TU Berlin\\
 \texttt{pokutta@zib.de}}
\date{\today}
\fi

\begin{document}
\ificmlbuild
\twocolumn[
  \icmltitle{Smoothness of the Feasible Set Does Not Accelerate Linear-Oracle Methods}
  \begin{icmlauthorlist}
    \icmlauthor{Sebastian Pokutta}{zib}
  \end{icmlauthorlist}
  \icmlaffiliation{zib}{Zuse Institute Berlin and Technische Universit\"at Berlin, Berlin, Germany}
  \icmlcorrespondingauthor{Sebastian Pokutta}{pokutta@zib.de}
  \icmlkeywords{Frank--Wolfe, linear minimization oracle, oracle complexity, lower bounds,
    smooth sets}
  \vskip 0.3in
]
\printAffiliationsAndNotice{}
\else
\maketitle
\fi

\begin{abstract}
Frank--Wolfe methods access the feasible region through a linear minimization oracle (LMO)
and reach primal gap $O(LD^2/T)$ after $T$ oracle calls on compact convex sets of diameter $D$,
the optimal order for such methods in general, and the accelerated rate $O(1/T^2)$ on
strongly convex sets. We show that smoothness of the set does not help beyond a trivial
reduction. For every $T\ge0$, $\beta>0$ and $D>2/\beta$, and every randomized method making at
most $T$ LMO calls, there is a $\beta$-smooth convex set of diameter at most $D$, on which
every nonzero linear function has a unique minimizer, such that with probability at least
$1/5$ the method returns an infeasible point or a point with primal gap at least
$L(D-2/\beta)^2/(64(T+1))$ for a known $L$-smooth and $L$-strongly convex quadratic.
A variant of Frank--Wolfe, which runs Frank--Wolfe on the set shrunk by $1/\beta$ and then
adds the ball back, reaches $2L(D-2/\beta)^2/(T+1)$ on every such set,
so the order $L(D-2/\beta)^2/T$ is optimal; this answers a question of Grimmer and Liu for
sets of constant smoothness. We add variants for methods with feasible outputs, for
deterministic methods and with a penalty for infeasible outputs, and upper bounds for
general smooth objectives and for sets with inner and outer tangent balls. In finite
dimension the hard set must depend on $T$; in general Hilbert spaces a single set suffices. One
smooth set with unique linear minimizers forces primal gap of order $1/t$ at every iteration
$t$ of every method that queries gradients only at points of the convex hull of its previous
answers, such as Frank--Wolfe; for every deterministic method with arbitrary queries,
one such set forces order $1/T$ after $T$ calls, for every horizon $T$; and for randomized
methods one set per method forces order $1/(T\log T(\log\log T)^2)$ at infinitely many
horizons. All theorems are formally verified in Lean~4.
\end{abstract}

\section{Introduction}\label{sec:intro}

We consider $\min_{x\in C}f(x)$ for a convex, $L$-smooth $f$ and a compact convex set $C$
that is accessible only through a \emph{linear minimization oracle} (LMO), which returns
for a direction $q$ a point $\LMO_C(q)\in\argmin_{z\in C}\ip qz$. Frank--Wolfe methods
\citep{frank1956algorithm,levitin1966constrained,CGFWSurvey2022} reach primal gap
$O(LD^2/T)$ after $T$ calls, $D=\diam C$, which is optimal over general compact convex sets
\citep{jaggi2013revisiting,lan2013complexity}. On strongly convex sets the rate improves to
$O(1/T^2)$ for strongly convex objectives \citep{garber2015faster}; \citet{KAP2020,KAP2021}
obtain faster rates on uniformly convex sets and \citet{kerdreux2022restarting} under H\"olderian
error bounds, and \citet{P2026FrankWolfe} shows $o(1/t)$ under local dual sharpness, a local form
of uniform convexity. Lower bounds of order $1/T^2$ on strongly convex sets are known for
Frank--Wolfe itself \citep{halbey2026lower}, for deterministic LMO methods
\citep{grimmer2026lower} and for randomized ones \citep{pokutta2026lower}.

A set is \emph{$\beta$-smooth} if its unit outer normals are $\beta$-Lipschitz along the
boundary (Definition~\ref{def:smooth}); smoothness is the dual notion of strong convexity, and methods
that access a smooth set through gauges, normals or membership queries exploit it
\citep{liu2023gauges,hazan2026linearly}. For LMO methods, \citet{grimmer2026lower} showed that
LMO-span methods are not accelerated on modestly smooth sets, $\beta=\Omega(1/\sqrt\eps)$, and left
open ``whether acceleration is possible when sets possess constant levels of smoothness'';
numerical performance estimation found no clear acceleration for small $T$
\citep{luner2026performance}.

\paragraph{Main result.} We answer the question of \citet{grimmer2026lower} negatively, and
not only for LMO-span methods: smoothness does not change the order $1/T$. Every
randomized LMO method with at most $T$ calls returns, with probability at least $1/5$, no
$\eps$-solution, that is, no feasible point with primal gap below
$\eps=L(D-2/\beta)^2/(64(T+1))$, on some $\beta$-smooth set of diameter at most $D$ with unique linear
minimizers, for a known $L$-smooth, $L$-strongly convex quadratic (Theorem~\ref{thm:main}), while
a variant of Frank--Wolfe, which runs Frank--Wolfe on the set shrunk by $1/\beta$, its inner
parallel body, and then adds the ball back (Definition~\ref{def:inner-fw}), reaches $2L(\diam C-2/\beta)^2/(T+1)$
for such quadratics on every
$\beta$-smooth set (Theorem~\ref{thm:upper}). Smoothness thus helps only by reducing the diameter from
$D$ to $D-2/\beta$ (Corollary~\ref{cor:optimal}).

The reason is simple. A $\beta$-smooth set is a convex set $K$ plus a ball of radius $1/\beta$,
and its LMO is an LMO of $K$ up to the shift $q/(\beta\norm q)$, which is known to the method
and costs no oracle call (Lemma~\ref{lem:rolling}).
Smoothness bounds the curvature from above, so faces may stay flat, and flat faces are where
LMO methods are slow. Our hard sets are rounded simplices whose objective is minimized next to
the barycenter of a flat face, and each LMO answer reveals at most one vertex of that face
(Figure~\ref{fig:hard-sets}). The main technical step removes all freedom of tie-breaking: perturbing
the vertices and thickening the simplex makes every linear minimizer unique, while each answer
still reveals at most one vertex (Section~\ref{sec:proof-main}). The resulting sets are smooth and
strictly convex. They are even strongly convex, but only with a modulus that tends to zero as
$T$ grows (Proposition~\ref{prop:strongly-convex}). At a fixed modulus Frank--Wolfe is faster on
uniformly convex sets \citep{KAP2020,KAP2021}, and sets that are smooth and strongly convex admit
a linear rate (Theorem~\ref{thm:tangent}; for gauge and normal oracles see
\citealp{liu2023gauges}). Tables~\ref{tab:lower} and~\ref{tab:upper} list the lower and the upper bounds.

\paragraph{Single sets.} Here and in the following, $T$ is the horizon, the number of oracle
calls after which a method returns its output, and $t$ counts the iterations of a running
method; a bound at every iteration $t$ holds along one run on one instance, an anytime bound.
The hard instance of Theorem~\ref{thm:main} depends on $T$ and has dimension $256(T+1)+1$. For
Frank--Wolfe this is unavoidable in finite dimension, where it converges at rate $o(1/t)$ on
every single set with unique linear minimizers once the gradient at the optimum is nonzero
(Proposition~\ref{prop:little-o}). In general Hilbert spaces, which may be infinite-dimensional, a
single set suffices (Section~\ref{sec:hilbert}). A single $1$-smooth set with unique linear
minimizers forces primal gap at least $31/(256(t+1))$ at every iteration $t$ of every
hull-based method started at a fixed point $x_0$ of the set (Definition~\ref{def:method}), which
queries gradients only at points of the convex
hull of $x_0$ and its previous answers, as the usual Frank--Wolfe variants do (Theorem~\ref{thm:hilbert}). A
method with arbitrary queries could instead query a direction whose answer reveals the
minimizer: on this set, the vertical direction (Remark~\ref{rem:exposure}). For arbitrary
queries the set must therefore depend on the method: for every deterministic method, a single
such set forces $1/(4(T+2))$ after $T$ calls, for every horizon $T$ (Theorem~\ref{thm:anytime}); and
for randomized methods, one distribution over such sets forces order
$1/(T\log T(\log\log T)^2)$ at every $T$, and for each randomized method one set, independent of its
randomness, forces this order at infinitely many $T$ (Theorem~\ref{thm:random}).

\paragraph{Formal verification.} All theorems, propositions, lemmas and corollaries are
kernel-checked theorems of Lean~4 \citep{moura2021lean} and Mathlib \citep{mathlib2020} in the
Lemmata Repository, entry \href{https://lemmata.iol.berlin/abs/2611.00021}{2611.00021}
\citep{lemmata2611}; only Remark~\ref{rem:two-term}, which relies
on classical first-order lower bounds, is not formalized. Formal results carry the mark
\leanverified.
The formal lower bounds of
Sections~\ref{sec:main}--\ref{sec:proof-main} live in $\R^N\times\R$: finitely many indexed
coordinates and one more coordinate, the height (Section~\ref{sec:setting}), the upper bounds hold
in every finite-dimensional inner product space, and Section~\ref{sec:hilbert} is formalized
in general real Hilbert spaces. We give full proofs in
Sections~\ref{sec:main}--\ref{sec:proof-main} and proof sketches for most results of
Sections~\ref{sec:beyond} and~\ref{sec:hilbert}.

\begin{widetable}[t]
\caption{Lower bounds. Each group first states the properties of the sets $C$; each row
then gives a bound that every method of the given kind satisfies on some such set, run with an
LMO of $C$. Unless stated otherwise the set depends on the method, but never on its randomness.
The objective is a known isotropic quadratic, $L$-smooth and $L$-strongly convex, with $L=1$ in
the last group, and $\delta=D-2/\beta$. ``No $\eps$-solution'' means that the output is infeasible
or has primal gap at least $\eps$. Randomized methods are measurable, in Theorem~\ref{thm:random}
jointly measurable (Definition~\ref{def:method}), and $c\,g(T)$ is of order
$1/(T\log T(\log\log T)^2)$.}
\label{tab:lower}
\centering\small
\setlength{\tabcolsep}{5pt}
\renewcommand{\arraystretch}{1.12}
\medskip
\begin{tabularx}{\textwidth}{@{}>{\raggedright\arraybackslash}p{4.7cm}Ylr@{}}
\toprule
Methods & After $T$ LMO calls & Probability & \\
\midrule
\multicolumn{4}{@{}p{\textwidth}@{}}{\textbf{Sets with unique linear minimizers} in
 $\R^{256(T+1)+1}$: $\beta$-smooth, $\diam C\le D$; strongly convex, but only with moduli at most
 $\beta/(10^8(T+1))$ (Proposition~\ref{prop:strongly-convex}). The set depends on $T$.}\\
\addlinespace[2pt]
randomized & no $\eps$-solution, $\eps=L\delta^2/(64(T+1))$ & ${\ge}\,1/5$ & Theorem~\ref{thm:main}\,(a)\\
randomized, feasible outputs & primal gap ${\ge}\,L\delta^2/(64(T+1))$ & ${\ge}\,3/4$ & Theorem~\ref{thm:main}\,(b)\\
\midrule
\multicolumn{4}{@{}p{\textwidth}@{}}{\textbf{Hard sets $C_I$ with the vertex LMO}
 (Definition~\ref{def:hard}) in $\R^{256(T+1)+1}$: $\beta$-smooth, $\diam C_I=D$, with a flat face
 and hence not strongly convex. The set depends on $T$.}\\
\addlinespace[2pt]
randomized & no $\eps$-solution, $\eps=L\delta^2/(32(T+1))$ & ${\ge}\,1/4$ & Theorem~\ref{thm:vertex}\\
randomized, feasible outputs & primal gap ${\ge}\,L\delta^2/(32(T+1))$ & ${\ge}\,3/4$ & Theorem~\ref{thm:variants}\,(a)\\
deterministic & no $\eps$-solution, $\eps=L\delta^2/(32(T+1))$, already in $\R^{64(T+1)+1}$ & --- & Theorem~\ref{thm:variants}\,(b)\\
randomized & penalized error ${\ge}\,L\delta^2/(256(T+1))$ & ${\ge}\,1/4$ & Theorem~\ref{thm:variants}\,(c)\\
\midrule
\multicolumn{4}{@{}p{\textwidth}@{}}{\textbf{Single sets in infinite-dimensional Hilbert spaces}:
 $1$-smooth, weakly compact, $\diam C\le15/4$, with unique linear minimizers. The set does not
 depend on $T$.}\\
\addlinespace[2pt]
hull-based, started at a fixed $x_0\in C$; one set for all methods & primal gap ${\ge}\,31/(256(t+1))$ at every iteration $t$ & --- & Theorem~\ref{thm:hilbert}\\
deterministic, arbitrary queries; one set per method & no $\eps$-solution for every $T$, $\eps=1/(4(T+2))$ & --- & Theorem~\ref{thm:anytime}\\
randomized, arbitrary queries; one distribution over sets for all methods & no $\eps$-solution, $\eps=c\,g(T)$, for every $T$ & ${\ge}\,1/5$ & Theorem~\ref{thm:random}\,(a)\\
randomized, arbitrary queries; one set per method & no $\eps$-solution, $\eps=c\,g(T)$, for infinitely many $T$ & ${\ge}\,1/5$ & Theorem~\ref{thm:random}\,(b)\\
\bottomrule
\end{tabularx}
\end{widetable}

\begin{widetable}[t]
\caption{Upper bounds on the primal gap. Each group first states the properties of the sets
$C$; each bound holds for every such set, every objective $f$ as stated and every LMO of $C$. Here
$D=\diam C$, $\delta=\diam C-2/\beta$, $K$ is the inner parallel body of $C$
(Lemma~\ref{lem:rolling}), ``smooth'' for $f$ means an $L$-Lipschitz gradient, and $c\in(0,1]$
depends only on $R/r$ (Theorem~\ref{thm:tangent}).}
\label{tab:upper}
\centering\small
\setlength{\tabcolsep}{5pt}
\renewcommand{\arraystretch}{1.12}
\medskip
\begin{tabularx}{\textwidth}{@{}>{\raggedright\arraybackslash}p{4.3cm}>{\raggedright\arraybackslash}p{3.7cm}Yr@{}}
\toprule
Method & Objective $f$ & Primal gap after $T$ LMO calls & \\
\midrule
\multicolumn{4}{@{}p{\textwidth}@{}}{\textbf{$\beta$-smooth sets}: compact and convex, in finite
 dimension.}\\
\addlinespace[2pt]
Frank--Wolfe on $K$ (Definition~\ref{def:inner-fw}) & isotropic quadratic $\frac L2\norm{\cdot-u}^2$ & ${\le}\,2L\delta^2/(T+1)$ for $T\ge2$ & Theorem~\ref{thm:upper}\\
conditional gradient sliding on $K\times B_{1/\beta}$ & convex, smooth & ${\le}\,540LD^2/T^2+810L\delta^2/T$ for $T\ge144$, where $T$ counts gradients and LMO calls together & Corollary~\ref{cor:two-term}\\
\midrule
\multicolumn{4}{@{}p{\textwidth}@{}}{\textbf{Sets with inner and outer tangent balls} of radii
 $r\le R$ at every boundary point: compact and convex, in finite dimension; these sets are smooth
 and strongly convex.}\\
\addlinespace[2pt]
projected gradient descent with Frank--Wolfe projections (Algorithm~\ref{alg:pgd-fw}) & $\mu$-strongly convex, smooth & ${\le}\,(1-c\mu/L)^{T-1}LD^2/2$, with $T$ gradients & Theorem~\ref{thm:tangent}\\
\midrule
\multicolumn{4}{@{}p{\textwidth}@{}}{\textbf{Norm-compact convex sets} in Hilbert spaces, in
 particular all compact convex sets in finite dimension.}\\
\addlinespace[2pt]
Frank--Wolfe in the sense of~\eqref{eq:fw-descent} & convex, smooth; $\nabla f(x^\star)\ne0$, and $x^\star$ is the only minimizer of $\ip{\nabla f(x^\star)}\cdot$ on $C$ & $o(1/t)$ at iteration $t$ & Proposition~\ref{prop:little-o}\\
\bottomrule
\end{tabularx}
\end{widetable}

\section{Setting}\label{sec:setting}

Throughout, $\ip\cdot\cdot$ is a real inner product with norm $\norm\cdot$, $B_r$ is the
closed ball of radius $r$ about the origin, $e_j$ are the unit vectors of $\R^N$,
$[N]=\{1,\dots,N\}$, and a \emph{convex body} is a compact convex set with nonempty interior.

\begin{definition}[Smooth sets]\label{def:smooth}
A vector $\nu$ with $\norm\nu=1$ is a \emph{unit outer normal} of a set $C$ at $z\in C$ if
$\ip\nu y\le\ip\nu z$ for all $y\in C$. The set $C$ is \emph{$\beta$-smooth} if
$\norm{\nu-\nu'}\le\beta\norm{z-z'}$ for all $z,z'\in C$, every unit outer normal $\nu$ at $z$ and
every unit outer normal $\nu'$ at $z'$.
\end{definition}

Each boundary point of a $\beta$-smooth set has exactly one unit outer normal, but a linear
function may still be minimized by a whole face. A $\beta$-smooth set is $\beta'$-smooth for
every $\beta'\ge\beta$; a set that is $\beta$-smooth but not $\beta'$-smooth for any $\beta'<\beta$ has
smoothness constant exactly $\beta$.

\begin{definition}[Linear minimization oracle]\label{def:lmo}
An \emph{LMO} of a set $C$ is a map $\LMO_C$ that returns for every direction $q$ a point
$\LMO_C(q)\in C$, the \emph{answer} to $q$, such that $\ip q{\LMO_C(q)}\le\ip qy$ for all $y\in C$.
\end{definition}

When several points minimize $\ip q\cdot$, an LMO selects one of them by a
\emph{tie-breaking rule}; for $q=0$ every point of $C$ is a minimizer.

\begin{definition}[Methods]\label{def:method}
We use the following classes of methods.
\begin{enumerate}
\item A \emph{deterministic method} requiring at most $T$ oracle calls either returns a point, or calls
the LMO with a direction $q$ and continues, depending on the answer, as a deterministic method requiring
at most $T-1$ calls; run with an LMO of a set $C$, it eventually returns a point of the
ambient space, its \emph{output}, which need not lie in $C$.
\item A \emph{randomized method} is a family of deterministic methods indexed by the outcome
$\omega$ of a random experiment with an arbitrary probability distribution. It is
\emph{measurable} if, for every measurable LMO, its output is a measurable function of $\omega$,
and \emph{jointly measurable} if its queries and its output are jointly Borel measurable functions
of $\omega$ and of the previous answers.
\item An \emph{anytime method} returns an output after every call: it chooses its $(T+1)$-st
query and its output after $T$ calls as functions of its first $T$ answers (and of $\omega$, if it
is randomized).
\item A \emph{hull-based method}, or \emph{Frank--Wolfe-type method}, for an objective $f$ and a start
$x_0$ calls in round $t=0,1,\dots$ the LMO with the gradient $\nabla f(z_t)$ at a point $z_t$ of the
convex hull of $x_0$ and the previous answers; its iterates and its output lie in this
hull. Frank--Wolfe with
the usual step sizes, short steps or exact line search, and its fully corrective, away-step and
pairwise variants are hull-based.
\end{enumerate}
\end{definition}

Directions and outputs may depend arbitrarily on previous answers. All lower bounds hold for
deterministic or randomized methods, which may query arbitrary directions, except
Theorem~\ref{thm:hilbert}, which concerns hull-based methods. The distinction matters for single sets
(Section~\ref{sec:hilbert}): for methods with arbitrary queries, a single set must depend on the method. In the lower bounds the
objective is known, so evaluating it or its gradient is free, and the hard instance may depend
on the method and on $T$ but not on $\omega$. The lower bounds of
Sections~\ref{sec:main}--\ref{sec:proof-main} use
\begin{equation}\label{eq:objective}
 f(x,y)=\frac L2\bigl(\norm x^2+y^2\bigr)+\frac{32L}{\beta}\,y,
 \qquad (x,y)\in\R^N\times\R ,
\end{equation}
which up to a constant is $\frac L2\norm{\cdot-c}^2$ with $c=(0,-32/\beta)$, so $L$-smooth and
$L$-strongly convex. For a point $(x,y)\in\R^N\times\R$ we call $x$ its \emph{horizontal part}
and the last coordinate $y$ its \emph{height}, and $(0,1)$ the \emph{vertical direction}; a point of
a set minimizes the vertical direction exactly when it has minimal height. A point $z$ is an
\emph{$\eps$-solution} on $C$ if $z\in C$ and $f(z)-\min_Cf<\eps$. We write $w$ for the output of a
method; the lower bounds show that $w$ is no $\eps$-solution, that is, $w\notin C$ or
$f(w)-\min_Cf\ge\eps$.

\begin{lemma}[Inner balls and the inner parallel body\leanverified]\label{lem:rolling}
Let $C$ be a nonempty compact convex $\beta$-smooth set in a finite-dimensional space
$E\ne\{0\}$, and $r=1/\beta$. (i)~For every unit outer normal $\nu$ of $C$ at $z\in C$,
$z-r\nu+B_r\subseteq C$. (ii)~The \emph{inner parallel body} $K=\{x:x+B_r\subseteq C\}$ satisfies
$C=K+B_r$ and $\diam C=\diam K+2r$. (iii)~For every $\LMO_C$ of $C$ and $q\ne0$, the point
$\LMO_C(q)+rq/\norm q$ lies in $K$ and minimizes $\ip q\cdot$ over $K$.
\end{lemma}
Parts (i) and~(ii) are classical: (i) is Blaschke's rolling theorem in the form of
\citet{walther1999generalization}, and (ii) is the characterization of balls as Minkowski summands
\citep[Section~3.2]{schneider2014convex}; extending a diametral segment of $K$ gives the diameter.
For~(iii), minimize over the two summands separately; on $B_r$ the minimizer $-rq/\norm q$ is
unique.

\section{Main results}\label{sec:main}

In the main result every nonzero direction has exactly one minimizer, so the LMO has no
choice to make, i.e., there are no ties to break.

\begin{definition}[The class $\cls(\beta,D,x_0)$]\label{def:class}
For $\beta,D>0$ and a point $x_0$, the class $\cls(\beta,D,x_0)$ consists of the $\beta$-smooth convex
bodies $C$ with $\diam C\le D$ and $x_0\in C$ on which $\ip q\cdot$ has exactly one minimizer
for every $q\ne0$; these are strictly convex. We write $\LMO_{C,x_0}$ for the LMO that returns this
minimizer for $q\ne0$ and $x_0$ for $q=0$.
\end{definition}

Every LMO of a set in $\cls(\beta,D,x_0)$ returns the same answer to every $q\ne0$, so no
tie-breaking rule is needed. For $q=0$ every point of the set is a minimizer, and the answer is
fixed ahead of time as $x_0$; otherwise an LMO could, for instance, return the minimizer of $f$
to $q=0$. With this convention the following bound holds for arbitrary LMOs.

\begin{theorem}[Lower bound for arbitrary LMOs\leanverified]\label{thm:main}
Let $T\ge0$, $L,\beta>0$ and $\beta D>2$. Put $N=256(T+1)$ and
$x_0=\bigl(0.99\,(D-2/\beta)\,e_1/\sqrt2,\,0\bigr)\in\R^N\times\R$, and let $\eps=L(D-2/\beta)^2/(64(T+1))$. Consider a
measurable randomized method requiring at most $T$ oracle calls.
\begin{enumerate}
\item Some $C\in\cls(\beta,D,x_0)$, depending on the method but not on $\omega$, makes the
 method, run with $\LMO_{C,x_0}$, return no $\eps$-solution with probability at least $1/5$.
\item If the method returns a point of $C$ when run with $\LMO_{C,x_0}$, for every $C\in\cls(\beta,D,x_0)$
 and every $\omega$, then some such $C$ gives $f(w)-\min_Cf\ge\eps$ with probability at least
 $3/4$, and expected primal gap at least $3\eps/4$.
\end{enumerate}
\end{theorem}

An infeasible output is never an $\eps$-solution. Frank--Wolfe methods return convex combinations of
oracle answers, so (b) applies to them. As $\diam C\le D$, the bounds also hold with $\diam C$
in place of $D$. The proof is outlined in Section~\ref{sec:proof-main} and given in
Appendix~\ref{app:main}. The matching upper bound runs
Frank--Wolfe on the inner parallel body~$K\subseteq C$ instead of on $C$. This loses nothing: by
Lemma~\ref{lem:rolling}, $C=K+B_r$, so minimizing $f$ over $C$ amounts to minimizing
$F(a)=\min_{b\in B_r}f(a+b)$ over $K$, and the LMO of $C$, shifted as in
Lemma~\ref{lem:rolling}\,(iii), is an LMO of $K$.

The method runs Frank--Wolfe on $K$ for $F$ and adds the optimal point of $B_r$ at the end.
For an isotropic quadratic $f=\frac L2\norm{\cdot-u}^2$ with known center $u$,
$F(a)=\frac L2(\norm{a-u}-r)_+^2$; here and below $(s)_+^2$ means $(\max\{s,0\})^2$, the positive part
taken before squaring, as $F(a)=0$ when $u\in a+B_r$. The gradient $\nabla F(a)=L(a+b_a-u)$, where
$a+b_a$ is the point of $a+B_r$ closest to $u$, is computed from $u$ in closed form; in this sense
the method needs no gradient evaluation, and its only oracle calls are LMO calls.

\begin{definition}[Frank--Wolfe on the inner parallel body]\label{def:inner-fw}
Let $f=\frac L2\norm{\cdot-u}^2$, $r=1/\beta$ and $\bar q\ne0$. Let $a_0=\LMO_C(\bar q)+r\bar q/\norm{\bar q}\in K$. For
$j=0,\dots,T-2$, call the LMO with $q_j=\nabla F(a_j)$, where $F(a)=\min_{b\in B_r}f(a+b)$, put
$s_j=\LMO_C(q_j)+rq_j/\norm{q_j}$ (stop if $q_j=0$) and $a_{j+1}=a_j+\frac2{j+2}(s_j-a_j)$. Return
$a+b$, where $a$ is the last iterate and $b$ the point of $B_r$ closest to $u-a$.
\end{definition}

We use the variant of Definition~\ref{def:inner-fw} instead of Frank--Wolfe on $C$ itself
because of the constant. Standard
Frank--Wolfe on $C$ with step sizes $2/(t+2)$ has primal gap at most $2L(\diam C)^2/{(T+2)}$ after
$T$ LMO calls \citep{jaggi2013revisiting}, with the full diameter of $C$, whereas the lower bound
of Theorem~\ref{thm:main} involves $D-2/\beta$. Our variant runs Frank--Wolfe only on the inner
parallel body $K$, whose diameter is $\diam C-2/\beta$, and handles the ball $B_r$ in closed form.
The same analysis then gives $2L(\diam C-2/\beta)^2/(T+1)$; the first LMO call provides the
starting point in $K$. Both bounds are of order $1/T$ and differ in the constant, by the factor
$(\diam C)^2/(\diam C-2/\beta)^2$, which is unbounded as $C$ approaches a ball of radius $1/\beta$.
This is not an artifact of the analysis: on a ball of radius $1/\beta$, $K$ is a single point and
our variant returns the minimizer, whereas standard Frank--Wolfe in
general does not reach it
after finitely many steps. So no bound proportional to $(\diam C-2/\beta)^2$ holds for standard
Frank--Wolfe. The variant of Definition~\ref{def:inner-fw}, in contrast, matches the lower
bound up to a constant for every $D>2/\beta$ (Corollary~\ref{cor:optimal}).

\begin{theorem}[Upper bound for isotropic quadratics\leanverified]\label{thm:upper}
Let $E$ be finite-dimensional, $T\ge2$, $L,\beta>0$, $u\in E$ and $f=\frac L2\norm{\cdot-u}^2$.
Frank--Wolfe on the inner parallel body requires at most $T$ oracle calls and no gradient
evaluation, and for every compact convex $\beta$-smooth $C\subseteq E$ with $\diam C\ge2/\beta$
and every LMO of $C$ it returns $w\in C$ with
\[
 f(w)-\min_Cf\;\le\;\frac{2L(\diam C-2/\beta)^2}{T+1}.
\]
\end{theorem}
\begin{proof}
By Lemma~\ref{lem:rolling}, $K$ has diameter $\delta=\diam C-2r$ and $\min_KF=\min_Cf$. The function
$F(a)=\frac L2(\norm{a-u}-r)_+^2$ is convex with gradient $L(a+b_a-u)$, where $a+b_a$ is the
point of $a+B_r$ closest to $u$; it is $L$-Lipschitz because $z\mapsto z-P_{B_r}(z)$ is
nonexpansive. The $s_j$ are exact linear minimizers on $K$, so the standard Frank--Wolfe recursion
\citep{jaggi2013revisiting}
$h_{j+1}\le(1-\gamma_j)h_j+\frac L2\gamma_j^2\delta^2$ with $\gamma_j=2/(j+2)$ gives
$h_j\le2L\delta^2/(j+2)$ for $h_j=F(a_j)-\min_KF$; the first call and $T-1$ steps give the bound.
\end{proof}

Since~\eqref{eq:objective} is an isotropic quadratic, the two theorems can be combined:

\begin{corollary}[Optimal order\leanverified]\label{cor:optimal}
Let $T\ge2$, $L,\beta>0$ and $D>2/\beta$, and let $N$ and $x_0$ be as in Theorem~\ref{thm:main}. On every
set of $\cls(\beta,D,x_0)$ and with every LMO, Frank--Wolfe on the inner parallel body returns
a point with primal gap at most $2L(D-2/\beta)^2/(T+1)$ for~\eqref{eq:objective}, whereas every
measurable randomized method with at most $T$ calls returns, on some set of $\cls(\beta,D,x_0)$,
no $\eps$-solution for $\eps=L(D-2/\beta)^2/(64(T+1))$ with probability at least $1/5$.
\end{corollary}

Hence $\Theta(L(D-2/\beta)^2/\eps)$ LMO calls are necessary and sufficient for primal gap $\eps$
on these sets in the worst case.

\section{Hard sets and the vertex LMO}\label{sec:vertex}

The hard sets are rounded simplices; an instance is an index set $I$, and the objective is
the same for all instances. All points of their flat bottom face minimize the vertical
direction, so there an LMO needs a tie-breaking rule. The vertex LMO breaks ties in the least
informative way: it returns a vertex, which reveals one index of $I$, whereas the barycenter of
the face would reveal all of $I$ (Figure~\ref{fig:hard-sets}a). Section~\ref{sec:proof-main} removes the
tie-breaking.

\begin{definition}[Hard sets and the vertex LMO]\label{def:hard}
Let $r=1/\beta$, $a=(D-2/\beta)/\sqrt2$, and for $I\subseteq[N]$ let $\Delta_I$ be the weight
vectors $\lambda\ge0$ in $\R^N$ supported on $I$ with $\sum_{i\in I}\lambda_i=1$. The
\emph{lifted simplex} is $S_I=\{(a\lambda,r):\lambda\in\Delta_I\}$ and the \emph{hard set} is
$C_I=S_I+B_r$. The \emph{vertex LMO} returns for $q=(q_x,q_y)\ne0$ the point $(ae_j,r)-rq/\norm q$,
where $j\in I$ has the smallest coordinate $q_{x,j}$ (the smallest such index in case of ties),
and for $q=0$ the point $(ae_1,0)$.
\end{definition}

\begin{widefigure}[t]
\centering
\includegraphics{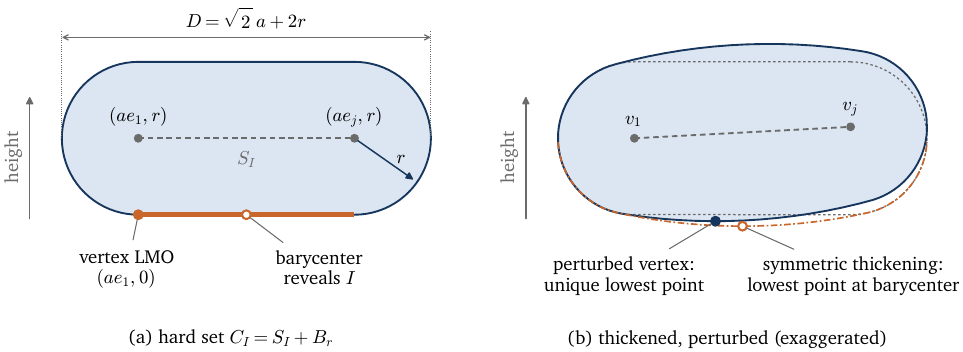}
\caption{Cross-sections through two vertices and the vertical direction, the direction of
the height (arrow). (a)~A hard set
(Definition~\ref{def:hard}): every point of the flat bottom face (orange) minimizes the vertical
direction; the vertex LMO returns a vertex and reveals one index of $I$, an answer at the
barycenter would reveal all of $I$. (b)~Thickening the simplex makes the lowest point unique,
but symmetric thickening puts it at the barycenter (Proposition~\ref{prop:symmetric}); with perturbed
vertices it again reveals a single vertex (Section~\ref{sec:proof-main}). Dotted: the set of~(a);
thickening and perturbation are exaggerated.}
\label{fig:hard-sets}
\end{widefigure}

For $q\ne0$ the vertex LMO minimizes $\ip q\cdot$ over each summand of $C_I$ separately, so it
returns a minimizer on $C_I$. For $q=0$ it returns the fixed point $(ae_1,0)$, which lies in $C_I$
when $1\in I$. So it is an LMO of $C_I$ whenever $1\in I$, and all index sets below contain $1$. A set $C$ is \emph{$(\alpha,2)$-uniformly convex} (strongly convex with modulus $\alpha$) if
$\gamma x+(1-\gamma)x'+\gamma(1-\gamma)\alpha\norm{x-x'}^2z\in C$ for all $x,x'\in C$,
$\gamma\in[0,1]$ and unit vectors $z$.

\begin{proposition}[Properties of the hard sets\leanverified]\label{prop:hard}
Let $1\in I$ and $|I|\ge2$. Then $C_I$ is a convex body of diameter exactly $D$; it is
$\beta$-smooth and not $\beta'$-smooth for any $\beta'<\beta$; for no $\alpha>0$ is it
$(\alpha,2)$-uniformly convex; and the vertex LMO is an LMO of $C_I$.
\end{proposition}
\begin{proof}
Write $z\in C_I$ as $z=v+b$ with $v\in S_I$, $\norm b\le r$. If $\nu$ is a unit outer normal at $z$,
comparing with $v+r\nu$ gives $b=r\nu$, and $v$ maximizes $\ip\nu\cdot$ on $S_I$. For $z'=v'+r\nu'$ this
gives $\ip{\nu-\nu'}{v-v'}\ge0$, hence $r\norm{\nu-\nu'}^2\le\ip{\nu-\nu'}{z-z'}\le\norm{\nu-\nu'}\norm{z-z'}$.
The points directly above and below a point of $S_I$ are $2r$ apart with opposite normals, so
no constant below $\beta$ works. The simplex has diameter $\sqrt2a$, so $\diam C_I\le\sqrt2a+2r=D$,
and extending a longest edge by $r$ at both ends attains it. The bottom points $(ae_i,0)$, $(ae_j,0)$ have their
midpoint at height $0$, the lowest height in $C_I$, so no ball about it lies in $C_I$.
\end{proof}

The faster rates of \citet{P2026FrankWolfe} hold under local dual sharpness: a set
$C$ with an $\LMO_C$ is \emph{locally dually sharp} around $M$ with constants $(A,p)$ if there is
$\rho>0$ with $A\norm g\norm{x-\LMO_C(g)}^p\le\ip g{x-\LMO_C(g)}$ for all $g$ and all $x\in C$ with
$\dist(x,M)<\rho$; Frank--Wolfe then converges faster than $1/t$. We now show that local dual
sharpness fails on the hard sets with the vertex LMO, so these faster rates do not contradict
the lower bounds of this section; for the sets of Theorem~\ref{thm:main} see
Proposition~\ref{prop:strongly-convex}.

\begin{proposition}[The minimizer\leanverified]\label{prop:optimum}
Let $|I|\ge2$, $c=(0,-32r)$, $o=(a\one_I/|I|,r)$ and $H=\norm{o-c}$. Then $z^\star=o+(r/H)(c-o)$
minimizes $f$ on $C_I$, and for the vertex LMO local dual sharpness around any $M\ni z^\star$
fails for all constants $A>0$ and exponents $p$.
\end{proposition}
\begin{proof}
Up to a constant $f=\frac L2\norm{\cdot-c}^2$. Every $v\in S_I$ has $\ip{v-o}{o-c}=0$, since
$o-c=(a\one_I/|I|,33r)$; so $z=v+b\in C_I$ satisfies
$H\norm{z-c}\ge\ip{z-c}{o-c}\ge H^2-rH$, i.e.\ $\norm{z-c}\ge\norm{z^\star-c}$. For $g=o-c$,
parallel to $\nabla f(z^\star)$, all $g_{x,i}$ with $i\in I$ are equal, the answer is
$s=(ae_j,r)-rg/H$ for one $j$, and $z^\star-s=o-(ae_j,r)\ne0$ while $\ip g{z^\star-s}=0$.
\end{proof}

The tie-breaking rule matters: for $q=(0,1)$ every point of the bottom face is a minimizer, and
an answer at its barycenter would reveal $I$ in one call. In contrast, the vertex LMO reveals little: a query
$q\ne0$ orders the indices by $q_{x,j}$, and the answer reveals only the first index of $I$ in
this order and that the indices before it are not in $I$; the answer to $q=0$ reveals nothing.

Fix a deterministic method, let $I$ be uniform over the index sets of size $d=8(T+1)$ that
contain $1$ (so that the vertex LMO is an LMO of $C_I$), and $N=256(T+1)$. After some answers, let $P$ be the indices known to lie in $I$,
$U$ the undecided ones, and $h=d-|P|$. Two index sets give the same answers exactly when they
agree on the decided indices; so, given the answers, $I$ is uniform over the sets $P\cup S$
with $S\subseteq U$, $|S|=h$, and an undecided index lies in $I$ with probability $h/|U|$. Each
call adds at most one index to $P$, so $|P|\le T+1$ at the end.

\begin{lemma}[Expected ratio of undecided indices\leanverified]\label{lem:posterior}
For a deterministic method requiring at most $T$ calls to the vertex LMO and $I$ uniform as
above, the final ratio satisfies $\E[h/|U|]={(d-1)/(N-1)}\le1/32$.
\end{lemma}
\begin{proof}
Each answer decides the indices before the first index of $I$ in its order, one at a time. A
decided index lies in $I$ with probability $h/|U|$, and the ratio becomes $(h-1)/(|U|-1)$;
otherwise it becomes $h/(|U|-1)$. The mean is $h/|U|$, so the expected ratio stays at its
initial value ${(d-1)/(N-1)}$.
\end{proof}

\begin{theorem}[Lower bound with the vertex LMO\leanverified]\label{thm:vertex}
Let $T\ge0$, $L,\beta>0$, $\beta D>2$ and $N=256(T+1)$. For every measurable randomized method
requiring at most $T$ calls to the vertex LMO there is an index set $I$ with $|I|=8(T+1)$ and
$1\in I$, depending on the method but not on $\omega$, such that, run with the vertex LMO of
$C_I$, the method returns no $\eps$-solution on $C_I$ for $\eps=L(D-2/\beta)^2/(32(T+1))$ with
probability at least $1/4$.
\end{theorem}
\begin{proof}
Write $k=T+1$, so $d=8k$. The barycenter $(a\one_I/d,0)$ of the bottom face gives
$\min_{C_I}f\le La^2/(2d)=La^2/(16k)$.

\emph{Step 1: little weight on $I$ is far from optimal.} For $w=(x,y)$ let $\theta_i=|x_i|/a$ if
$\norm x^2\le a^2/(4k)$ and $\theta=0$ otherwise, so $k\norm\theta^2\le1/4$. Suppose $w\in C_I$
and $\sum_{i\in I}\theta_i\le3/4$. Points of $C_I$ have $y\ge0$, so if $\norm x^2>a^2/(4k)$, then
$f(w)\ge\frac L2\norm x^2\ge La^2/(8k)$. Otherwise some $\lambda\in\Delta_I$ has
$\norm{x-a\lambda}^2+(y-r)^2\le r^2$, so $\norm{x-a\lambda}^2\le2ry$; as $\sum_{i\in I}x_i\le3a/4$,
$\sum_{i\in I}(a\lambda_i-x_i)\ge a/4$, and Cauchy--Schwarz over $I$ gives
$a^2/16\le d\norm{x-a\lambda}^2\le2dry$, so $f(w)\ge32Lry\ge La^2/(8k)$. Either way
$f(w)-\min_{C_I}f\ge La^2/(16k)$.

\emph{Step 2: little weight on $I$ is likely.} Fix $\omega$ and let $I$ be uniform as above.
Given the final answers, which determine the output, $i\in I$ with probability $p_i$: $p_i=1$
on $P$, $p_i=h/|U|$ on $U$, and $0$ otherwise. By $p_i\theta_i\le p_i^2/(4k)+k\theta_i^2$, the
expected weight on $I$ is at most
\[
 \frac{\E\bigl[|P|+h\cdot h/|U|\bigr]}{4k}+k\norm\theta^2
 \;\le\;\frac{k+d\,\E[h/|U|]}{4k}+\frac14\;\le\;\frac5{16}+\frac14=\frac9{16}
\]
by Lemma~\ref{lem:posterior}. By Markov's inequality the weight on $I$ is at most $3/4$ with
probability at least $1/4$, and then the output is no $\eps$-solution by Step~1.

\emph{Step 3.} Averaging over $\omega$, the probability of no $\eps$-solution, averaged over $I$, is at least
$1/4$; some $I$ attains the average, and $La^2/(16k)=L(D-2/\beta)^2/(32(T+1))$.
\end{proof}

Steps~1 and~2 only use that each call reveals the first index of $I$ in an order chosen by the
method. Figure~\ref{fig:rates} shows Frank--Wolfe on the inner parallel body between the upper bound
of Theorem~\ref{thm:upper} and the lower bound of Theorem~\ref{thm:vertex}. The
same family gives three further lower bounds, proved in Appendix~\ref{app:variants}.

\begin{widefigure}[t]
\centering
\includegraphics{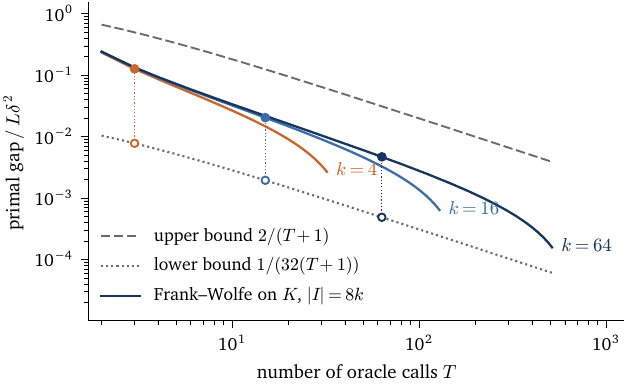}
\caption{Frank--Wolfe on the inner parallel body (Definition~\ref{def:inner-fw}) with the vertex LMO on
hard sets with $|I|=8k$, $N=256k$, $a=r=1$, run until all indices of $I$ are found. Dashed: the
upper bound of Theorem~\ref{thm:upper}; dotted: the lower bound $\eps$ of Theorem~\ref{thm:vertex}, which
holds for instances depending on $T$. At the horizon $T=k-1$ of each instance the primal gap (filled
markers) lies above this bound (open markers).}
\label{fig:rates}
\end{widefigure}

\begin{theorem}[Variants with the vertex LMO\leanverified]\label{thm:variants}
In the setting of Theorem~\ref{thm:vertex}, with index sets $I\ni1$ of size $8(T+1)$ depending on the
method but not on $\omega$, and $\eps=L(D-2/\beta)^2/(32(T+1))$:
\begin{enumerate}
\item If the method returns a point of $C_I$ for every index set $I\ni1$ and every $\omega$, then
 some $I$ gives $f(w)-\min_{C_I}f\ge\eps$ with probability at least $3/4$, and expected primal
 gap at least $3\eps/4$.
\item If $N=64(T+1)$, every deterministic method returns no $\eps$-solution on some $C_I$.
\item For $\rho\ge L(2D+32/\beta)$ and the penalized error $\err_I(w)=f(w)-\min_{C_I}f+\rho\dist(w,C_I)$,
 some $I$ gives $\err_I(w)\ge\eps/8$ with probability at least $1/4$ and $\E\err_I(w)\ge\eps/32$,
 without any feasibility assumption.
\end{enumerate}
\end{theorem}

The penalized error is the measure suggested by \citet{grimmer2026lower}. The argument
for~(b) does not extend to randomized methods, whose revealed indices depend on $\omega$; this is
why the randomized bounds use the larger dimension. Section~\ref{sec:proof-main} removes the tie-breaking rule and allows
arbitrary LMOs.

\section{Removing tie-breaking: arbitrary LMOs}\label{sec:proof-main}

To dispense with tie-breaking, every nonzero direction needs a unique minimizer. Thickening
the simplex achieves this: every point $\sum_{j\in I}\lambda_j(Ae_j,1)$ of the simplex is replaced
by a ball of radius $\sigma(1-\norm\lambda^2)$, which vanishes at the vertices and is largest at the
barycenter; the union of these balls is strictly convex, and adding $B_1$ makes it smooth
(Figure~\ref{fig:thickening}). But symmetric thickening reveals $I$ (Figure~\ref{fig:hard-sets}b):
recall that the height is the last coordinate, so the lowest points of a set are its points of
minimal height, the minimizers of the vertical direction.

\begin{widefigure}[t]
\centering
\includegraphics{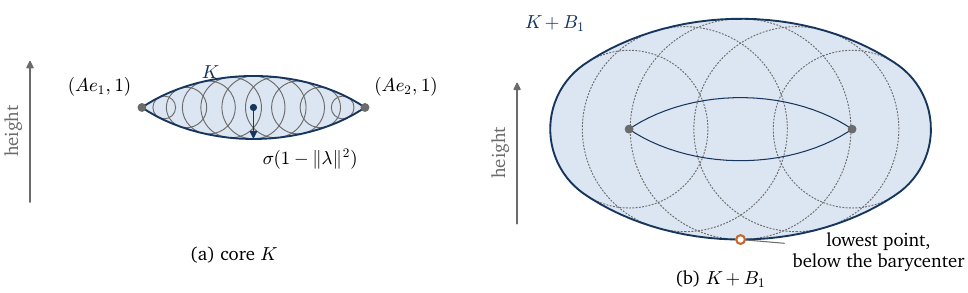}
\caption{The thickened sets of Section~\ref{sec:proof-main} in a cross-section through two vertices; the
arrow is the vertical direction, and $\sigma$ is exaggerated. (a)~The core $K$ is the union of the
balls of radius $\sigma(1-\norm\lambda^2)$ about the points of the simplex: the radius vanishes at
the vertices and is largest at the barycenter. (b)~$K+B_1$ is the union of the balls of radius
$1+\sigma(1-\norm\lambda^2)$ about the same points; with symmetric thickening its lowest point lies
below the barycenter (Proposition~\ref{prop:symmetric}), so one call in the vertical direction reveals~$I$.}
\label{fig:thickening}
\end{widefigure}

\begin{proposition}[Symmetric thickening reveals the index set\leanverified]\label{prop:symmetric}
Let $A,\sigma>0$ and let $K$ be the union of the balls of radius $\sigma(1-\norm\lambda^2)$ about
$\sum_{j\in I}\lambda_j(Ae_j,1)$, $\lambda\in\Delta_I$. Every point of minimal height of $K+B_1$ has horizontal
coordinates $A/|I|$ on $I$ and $0$ elsewhere.
\end{proposition}
\begin{proof}
The vertices have height $1$, so the lowest point over $\lambda$ has height
$-\sigma(1-\norm\lambda^2)$, minimized only at $\lambda=\one_I/|I|$.
\end{proof}

So one LMO call in the vertical direction reveals $I$, and we need a more elaborate
construction: the vertices must also be perturbed; then the best vertex is unique and, as with the vertex
LMO, reveals one index.

We take $r=L=1$, since scaling by $r$ multiplies distances by $r$ and primal gaps by $Lr^2$,
and put $k=T+1$, $N=256k$, $n=N+1$, $A_0=(\beta D-2)/\sqrt2$ and $A=0.99A_0$. A \emph{table}
$Z$ perturbs each vertex $(Ae_j,1)$, $j\ne1$, by a random vector $Z_j$ with coordinates on a
fine grid and $\norm{Z_j}\le\eta=A^2/(10^6k(A+36)^2)$, and with $v_j=(Ae_j,1)+Z_j$ and a tiny
$\sigma>0$ (both fixed in Appendix~\ref{app:main}) we consider
\[
 C_{I,Z}=\bigcup_{\lambda\in\Delta_I}\Bigl(\sum_{j\in I}\lambda_jv_j+\bigl(1+\sigma(1-\norm\lambda^2)\bigr)B_1\Bigr).
\]
These are $1$-smooth convex bodies of diameter at most $\beta D$ with unique linear minimizers,
since minimizing $\ip q\cdot$ amounts to minimizing the strictly convex function
$\sum_{j}\lambda_j\ip q{v_j}+\sigma\norm\lambda^2$ over $\Delta_I$.

\begin{proofsketch}[Proof sketch of Theorem~\ref{thm:main}]
The proof, given in Appendix~\ref{app:main}, has three steps.
(1)~\emph{One index per answer.} If the best index $j^\star$ of $I$ in the direction $q$ beats
every not yet revealed index of $I$ by more than $2\sigma$, the minimizing weights vanish outside the
revealed indices and $j^\star$. The method is then simulated by a method that learns the first
index of $I$ in an order of its choice, to which Steps~1--2 of the proof of Theorem~\ref{thm:vertex}
apply. (2)~\emph{Ambiguity is rare.} A comparison is ambiguous for at most one grid value of one
coordinate of $Z$, and a union bound over the decision tree leaves at most a fraction $1/100$ of
bad tables. (3)~\emph{Transfer.} The perturbed and unperturbed sets are within Hausdorff
distance $2\eta$, and a contraction toward $v_1$, known to the method, moves outputs between them at a cost of
at most $4A^2/(10^6k)$ in primal gap. So on a fraction $\frac14\cdot\frac{99}{100}\ge\frac15$ of the
instances the method returns no $\eps$-solution for $\eps=A_0^2/(32k)$, which is $L(D-2/\beta)^2/(64(T+1))$
after scaling; for~(b), the argument of Theorem~\ref{thm:variants}\,(a) gives
$\frac{97}{128}\cdot\frac{99}{100}\ge\frac34$.
\end{proofsketch}

The sets $C_{I,Z}$ have unique linear minimizers, so they are strictly convex. The thickening
even makes them strongly convex, but only barely.

\begin{proposition}[Strong convexity of the sets of Theorem~\ref{thm:main}\leanverified]\label{prop:strongly-convex}
In Theorem~\ref{thm:main}\,(a) and~(b) the set $C$ can be chosen such that, in addition, it is
$(\alpha,2)$-uniformly convex for some $\alpha>0$, and every $\alpha$ for which it is
$(\alpha,2)$-uniformly convex satisfies $\alpha\le\beta/(10^8(T+1))$.
\end{proposition}
\begin{proof}
The sets are the sets $rC_{I,Z}$ of the proof of Theorem~\ref{thm:main}. We show both claims
for $C_{I,Z}$; scaling by $r=1/\beta$ multiplies the moduli by $\beta$. Write
$c(\lambda)=\sum_{j\in I}\lambda_jv_j$ and $R(\lambda)=1+\sigma(1-\norm\lambda^2)$, so that $C_{I,Z}$ is
the union of the balls of radius $R(\lambda)$ about $c(\lambda)$, $\lambda\in\Delta_I$; the parameters of
Appendix~\ref{app:main} satisfy $\sigma\le\eta/8$ and $\eta\le10^{-6}$.

\emph{A positive modulus.} Let $x=c(\lambda)+p$ and $x'=c(\lambda')+p'$ be points of $C_{I,Z}$, with
$\norm p\le R(\lambda)$ and $\norm{p'}\le R(\lambda')$, let $\gamma\in[0,1]$, and put $\tau=\gamma(1-\gamma)$
and $\bar R=\gamma R(\lambda)+(1-\gamma)R(\lambda')$. The ball of the mixed weights
$\gamma\lambda+(1-\gamma)\lambda'$ has radius $\bar R+\sigma\tau\norm{\lambda-\lambda'}^2$, and
\[
 \norm{\gamma p+(1-\gamma)p'}^2=\gamma\norm p^2+(1-\gamma)\norm{p'}^2-\tau\norm{p-p'}^2
 \;\le\;\bar R^2+\tau\bigl((R(\lambda)-R(\lambda'))^2-\norm{p-p'}^2\bigr),
\]
where $|R(\lambda)-R(\lambda')|=\sigma\,\bigl|\norm\lambda^2-\norm{\lambda'}^2\bigr|\le2\sigma\norm{\lambda-\lambda'}$.
As $1\le\bar R\le5/4$ and $\sigma\le1/4$, this gives
$\norm{\gamma p+(1-\gamma)p'}\le\bar R+\tau\bigl(\frac\sigma2\norm{\lambda-\lambda'}^2-\frac25\norm{p-p'}^2\bigr)$.
So the ball of the mixed weights contains the ball of radius
$\tau\bigl(\frac\sigma2\norm{\lambda-\lambda'}^2+\frac25\norm{p-p'}^2\bigr)$ about $\gamma x+(1-\gamma)x'$. As
$\norm{x-x'}^2\le2V\norm{\lambda-\lambda'}^2+2\norm{p-p'}^2$ for $V=\sum_{j\in I}\norm{v_j}^2$, this radius is at
least $\tau\alpha\norm{x-x'}^2$ for $\alpha=\min\{1/5,\sigma/(4(V+1))\}$, and $C_{I,Z}$ is
$(\alpha,2)$-uniformly convex.

\emph{Every modulus is small.} The vertices $v_j$ have heights between $1-\eta$ and $1+\eta$. So every
point of $C_{I,Z}$ has height at least $-(\eta+\sigma)$. For two different vertices $v_i$ and $v_j$,
the lowest points $v_i-(0,1)$ and $v_j-(0,1)$ of the unit balls about them lie in $C_{I,Z}$ and
have height at most $\eta$. If $C_{I,Z}$ is $(\alpha,2)$-uniformly convex, it also contains the
point obtained from the midpoint of these two points by lowering its height by
$\alpha\norm{v_i-v_j}^2/4$; the height of this point is at most $\eta-\alpha\norm{v_i-v_j}^2/4$ and at
least $-(\eta+\sigma)$, hence $\alpha\norm{v_i-v_j}^2/4\le2\eta+\sigma\le3\eta$. As
$\norm{v_i-v_j}\ge A-2\eta\ge\frac{1999}{2000}A$ and $\eta=A^2/(10^6k(A+36)^2)$, this gives
$\alpha\le12.02/(10^6k(A+36)^2)\le1/(10^8k)$.
\end{proof}

\section{Beyond isotropic quadratics}\label{sec:beyond}

For a general smooth convex objective, conditional gradient sliding \citep{lan2016conditional}
applies to $K\times B_r$, with one gradient of $f$ per outer iteration and the LMO of $K$ from
Lemma~\ref{lem:rolling}\,(iii) in the inner Frank--Wolfe loops. The method and its analysis are those of
\citet{lan2016conditional}; the reduction to $K\times B_r$ confines the LMO calls to the factor $K$,
of diameter $\delta$.

\begin{theorem}[Smooth convex objectives, after \citet{lan2016conditional}\leanverified]\label{thm:sliding}
Let $E$ be finite-dimensional, $f$ convex and differentiable with $L$-Lipschitz gradient, $C$
compact, convex and $\beta$-smooth with $\diam C\ge2/\beta$, $\delta=\diam C-2/\beta$ and
$\Lambda=\delta^2+4/\beta^2$. For $m\ge1$, conditional gradient sliding on $K\times B_{1/\beta}$
evaluates $m$ gradients, makes at most $1+m+9\delta^2m(m+1)/\Lambda$ LMO calls and returns a
point of $C$ with primal gap at most $15L\Lambda/((m+1)(m+2))$.
\end{theorem}
\begin{proofsketch}
We follow the analysis of \citet{lan2016conditional}. The function $g(a,b)=f(a+b)$ is $2L$-smooth on $K\times B_r$, of squared diameter $\Lambda$.
With $\gamma_j=3/(j+2)$, $\tau_j=6L/(j+1)$ and tolerances $\eta_j=2L\Lambda/(j(j+1))$, the
three-point inequality gives
$h_j\le(1-\gamma_j)h_{j-1}+\gamma_j\eta_j+\frac{\gamma_j\tau_j}2(d_{j-1}-d_j)$ for the primal gaps
$h_j$ and the squared distances $d_j$ of the proximal iterates from a minimizer, which telescopes
with weights $j(j+1)(j+2)$. The proximal subproblem is solved exactly on $B_r$ and on $K$ by
Frank--Wolfe with exact line search until the Frank--Wolfe gap is at most $\eta_j$, within
$1+6\tau_j\delta^2/\eta_j=1+18j\delta^2/\Lambda$ LMO calls.
\end{proofsketch}

\begin{corollary}[Gradients and LMO calls together\leanverified]\label{cor:two-term}
In the setting of Theorem~\ref{thm:sliding} with $D=\diam C$, for every integer $T\ge144$ the sliding
method with $m=\lceil(15L\Lambda/\eps_T)^{1/2}\rceil$ makes at most $T$ gradient evaluations and
LMO calls in total and returns a point of $C$ with primal gap at most
$\eps_T=540LD^2/T^2+810L\delta^2/T$.
\end{corollary}

\begin{remark}\label{rem:two-term}
The second term matches Theorem~\ref{thm:main} up to constants. A lower bound of order $LD^2/T^2$ for
methods that pay for gradients follows from classical first-order lower bounds
\citep{nemirovski1983problem,drori2017exact} on a ball, whose LMO carries no information; this
is not formalized in Lean.
\end{remark}

\begin{algorithm2e}[t]
\caption{Projected gradient descent with Frank--Wolfe projections}\label{alg:pgd-fw}
\KwIn{$x_0\in C$, an $\LMO_C$ of $C$, the gradient $\nabla f$, $L$, $r$, and $m\ge0$}
$q\gets\nabla f(x_0)$\;
\lIf{$q=0$}{\Return{$x_0$}}
$x\gets\LMO_C(q)$; $a\gets x+rq/\norm q$\Comment*{$a\in K$ and $x\in a+B_r$}
\For{$j=1,\dots,m$}{
  $u\gets x-\nabla f(x)/L$\;
  \If{$\norm{u-a}>r$}{
    $s\gets\LMO_C(a-u)+r(a-u)/\norm{a-u}$\Comment*{answer on $K$, Lemma~\ref{lem:rolling}\,(iii)}
    $\gamma\gets\min\{1,\max\{0,\ip{a-u}{a-s}/\norm{s-a}^2\}\}$ ($\gamma=0$ if $s=a$)\;
    $a\gets a+\gamma(s-a)$\Comment*{exact line search for $\frac12\norm{\cdot-u}^2$}
  }
  $x\gets$ the point of $a+B_r$ closest to $u$\;
}
\Return{$x$}\;
\end{algorithm2e}

If $C$ also has an outer tangent ball at every boundary point and $f$ is strongly convex, a linear rate holds. Such sets
are smooth and strongly convex. Linear rates on them are known for methods with gauge and
normal oracles \citep{liu2023gauges}, and on smooth sets for methods with membership queries
\citep{hazan2026linearly}; with an LMO, Frank--Wolfe converges linearly on strongly convex sets when
the gradient is bounded away from zero \citep{levitin1966constrained}, and at rate $O(1/T^2)$ in
general \citep{garber2015faster}. With both tangent balls, one gradient and one LMO call per
iteration give a linear rate in general.
\emph{Projected gradient descent with Frank--Wolfe projections} (Algorithm~\ref{alg:pgd-fw}) keeps
$a\in K$ and $x\in a+B_r$. Each iteration takes the gradient step $u=x-\nabla f(x)/L$ and replaces
the projection of $u$ onto $C$ by a single Frank--Wolfe step with exact line search on $K$ for
$\frac12\norm{\cdot-u}^2$, followed by the point of $a+B_r$ closest to $u$. Conditional gradient
sliding \citep{lan2016conditional} also computes projections by Frank--Wolfe steps, but runs
several of them per projection, inside an accelerated outer loop. Here a single step suffices:
the outer tangent balls make it reduce the excess $\norm{a-u}^2-\min_{v\in K}\norm{v-u}^2$ by a
constant factor.

\begin{theorem}[Inner and outer tangent balls\leanverified]\label{thm:tangent}
Let $E\ne\{0\}$ be finite-dimensional, $C\subseteq E$ compact and convex, $0<r\le R$, and
suppose every $z\in C$ with unit outer normal $\nu$ satisfies $z-r\nu+B_r\subseteq C\subseteq
z-R\nu+B_R$. Let $f$ be differentiable and $\mu$-strongly convex, $\mu>0$, with $L$-Lipschitz
gradient, and $c=\frac12\min\{1,r/(2(R-r))\}$ if $R>r$, $c=1$ if $R=r$. From any $x_0\in C$,
projected gradient descent with Frank--Wolfe projections uses in $m$ iterations at most $m+1$ gradients
and $m+1$ LMO calls and returns $x_m\in C$ with
$f(x_m)-\min_Cf\le(1-c\mu/L)^m\,L(\diam C)^2/2$. If $R=r$, a variant uses one LMO call in total.
\end{theorem}
\begin{proofsketch}
For the answer $s\in K$ to $q=a-u$, the outer ball gives $2(R-r)G\ge\norm q\norm{s-a}^2$ for the
Frank--Wolfe gap $G=\ip q{a-s}$ of $\frac12\norm{\cdot-u}^2$ at $a$, so the point $a'$ found by
the line search satisfies $\norm{a'-u}^2\le(1-c)\norm{a-u}^2+c\norm{v-u}^2$ for every $v\in K$. The convex nondecreasing map
$t\mapsto\frac12(\max\{0,\sqrt{2t}-r\})^2$ turns this into $Q(x')\le(1-c)Q(x)+cQ(z)$ for
$z\in C$, the new iterate $x'$ and $Q=\frac12\norm{\cdot-u}^2$, and strong convexity and smoothness give the factor
$1-c\mu/L$.
\end{proofsketch}

\section{Single sets in general Hilbert spaces}\label{sec:hilbert}

The instances of Theorem~\ref{thm:main} depend on the method and on $T$. In general Hilbert
spaces, which may be infinite-dimensional, a single set suffices. For hull-based methods
(Definition~\ref{def:method}) one set works for every method and
every iteration $t$ (Theorem~\ref{thm:hilbert}); for deterministic methods with arbitrary queries one
set per method works for every horizon (Theorem~\ref{thm:anytime}); for randomized methods one
distribution over sets works for every horizon, and one set per method for infinitely many
horizons, up to logarithmic factors (Theorem~\ref{thm:random}). All these sets are weakly compact,
$1$-smooth, and have unique linear minimizers. A set is \emph{norm compact} if it is compact in
the norm topology, and \emph{weakly compact} if it is compact in the weak topology. In finite
dimension both are the usual compactness; in infinite dimension every closed bounded convex set
is weakly compact but need not be norm compact, the closed unit ball for instance is not. As in
the introduction, $t$ counts the
iterations of a running method, and a bound at every $t$ is an anytime bound along one run; $T$
is the horizon, the number of calls after which a method returns its output.

\subsection{Hull-based methods}

Here \emph{Frank--Wolfe} is any iteration with $s_t=\LMO_C(\nabla f(x_t))$ and $x_{t+1}\in C$ such
that, with $\gamma_t=2/(t+2)$,
\begin{equation}\label{eq:fw-descent}
 f(x_{t+1})\le f(x_t)-\gamma_t\ip{\nabla f(x_t)}{x_t-s_t}+\frac L2\gamma_t^2\norm{s_t-x_t}^2;
\end{equation}
by $L$-smoothness this covers the step sizes $2/(t+2)$, short steps, exact line search and fully
corrective Frank--Wolfe. On norm-compact sets with unique linear minimizers, in particular in
finite dimension, a single
set cannot force the order $1/t$ on these methods when the gradient at the minimizer is
nonzero. The lower bound of
Theorem~\ref{thm:hilbert} below holds for all hull-based methods, which need not satisfy
\eqref{eq:fw-descent}; both classes contain the usual Frank--Wolfe variants.

\begin{proposition}[Norm-compact sets\leanverified]\label{prop:little-o}
Let $C$ be a norm-compact convex set in a real Hilbert space, for instance a compact convex set
in $\R^n$, let $f$ be convex and differentiable with $L$-Lipschitz gradient, and let
$x^\star\in\argmin_Cf$ with $g^\star=\nabla f(x^\star)\ne0$ such that $x^\star$ is the only minimizer
of $\ip{g^\star}\cdot$ on $C$. Then Frank--Wolfe with any LMO of $C$ satisfies
$t\,(f(x_t)-f(x^\star))\to0$.
\end{proposition}
\begin{proof}
The usual Frank--Wolfe analysis gives $h_t=f(x_t)-f(x^\star)\to0$. Convexity and optimality give
$0\le\ip{g^\star}{x_t-x^\star}\le h_t\to0$. By compactness and uniqueness, every sequence $y_t\in C$
with $\ip{g^\star}{y_t-x^\star}\to0$ converges to $x^\star$; so $x_t\to x^\star$, $\nabla f(x_t)\to
g^\star$, and $s_t\to x^\star$ since $\ip{g^\star}{s_t-x^\star}\le\norm{g^\star-\nabla f(x_t)}\diam C$.
Thus $\ell_t=\norm{s_t-x_t}\to0$. With $\ip{\nabla f(x_t)}{x_t-s_t}\ge h_t$, \eqref{eq:fw-descent}
gives $h_{t+1}\le(1-\gamma_t)h_t+\frac L2\gamma_t^2\ell_t^2$, i.e.\ $(t+1)(t+2)h_{t+1}\le
t(t+1)h_t+2L\ell_t^2$, and summing, $(t+1)h_{t+1}\le\frac{2L}{t+2}\sum_{s\le t}\ell_s^2\to\nobreak0$.
\end{proof}

\citet{P2026FrankWolfe} proves $o(1/t)$ under local dual sharpness, a power-type condition,
and shows with the example of a planar body whose boundary near a support point is the graph
of $e^{-1/x^2}$ that strict convexity alone does not imply it. Proposition~\ref{prop:little-o} assumes only uniqueness, gives
no rate, and applies to such bodies; local dual sharpness, conversely, allows non-unique
minimizers and, together with a local H\"olderian error bound, yields explicit rates.

The hard sets of Theorem~\ref{thm:main} satisfy the assumptions of
Proposition~\ref{prop:little-o}: they have unique linear minimizers, and~\eqref{eq:objective} has
nonzero gradient at its minimizer, as its center $c$ lies outside of them. In infinite dimension a unique minimizer of
$\ip{g^\star}\cdot$ need not attract the sequences along which $\ip{g^\star}{\cdot-x^\star}$
vanishes, and a single set can force the order $1/t$. Let $\tau_i=2^{-i-5}$ for $i\ge1$ and $\sigma=1/100$, and let
$d_1,d_2,\dots$ be vectors of a real Hilbert space $\mathcal H_0$ with Gram matrix
$\ip{d_i}{d_j}=\mathbf 1_{i=j}-\tau_{\max(i,j)}$, for instance $d_i=(\Id-\Gamma)^{1/2}e_i$ in $\ell_2$,
where $\Gamma_{ij}=\tau_{\max(i,j)}$. Since $\Gamma=\sum_j(\tau_j-\tau_{j+1})u_ju_j^\top$ with
$u_j=\sum_{i\le j}e_i$, $\Gamma\succeq0$ and $\norm\Gamma\le\sum_j(\tau_j-\tau_{j+1})j=\sum_j\tau_j=1/32$;
hence $S\xi=\sum_i\xi_id_i$ satisfies $\frac{31}{32}\norm\xi^2\le\norm{S\xi}^2\le\norm\xi^2$ for
finitely supported $\xi$. Let $\Delta=\{\lambda\in\ell_2:\lambda\ge0,\ \sum_i\lambda_i\le1\}$ and
\[
 c(\lambda)=\Bigl(S\lambda,\ \sum_i\tau_i^2\lambda_i\Bigr),\quad
 R(\lambda)=1+\sigma\sum_i\tau_i^2(\lambda_i-\lambda_i^2),\quad
 C=\bigcup_{\lambda\in\Delta}\bigl(c(\lambda)+R(\lambda)B_1\bigr)\subseteq\mathcal H_0\times\R,
\]
and $f(x,y)=\frac12(\norm x^2+(y+4)^2)$, which is $1$-smooth and $1$-strongly convex.

\begin{lemma}[The set\leanverified]\label{lem:hilbert-set}
(i)~$C$ is convex and weakly compact, contains $B_1$, has diameter at most $15/4$, is
$1$-smooth, and every nonzero linear function has exactly one minimizer on $C$.
(ii)~$z^\star=(0,-1)$ minimizes $f$ on $C$, $\nabla f(z^\star)=(0,3)$, and
$f(x,y)-f(z^\star)\ge\frac12\norm x^2$ on $C$.
\end{lemma}
\begin{proofsketch}
As $c$ is affine and $R$ concave, convex combinations of points of two balls lie in the
ball of the mixed parameter; with $R-1\ge0$ in place of $R$ this gives $C=K+B_1$ for a convex
$K$, hence $1$-smoothness as in the proof of Proposition~\ref{prop:hard}. Weak compactness follows since
$\Delta$ is weakly compact, $c$ is weakly continuous and $\lambda\mapsto\sum_i\tau_i^2\lambda_i^2$ is a
compact quadratic form. For $q\ne0$, minimizing $\ip q\cdot$ amounts to maximizing a strictly
concave function of $\lambda$ over $\Delta$. Points of $C$ have height at least
$-1+(1-\sigma)\sum_i\tau_i^2\lambda_i\ge-1$, which gives~(ii).
\end{proofsketch}

\begin{theorem}[One set for all iterations\leanverified]\label{thm:hilbert}
Let $q_0=(d_1/8,-1)$ and let $x_0$ be the maximizer of $\ip{q_0}\cdot$ on $C$.
\begin{enumerate}
\item Consider any hull-based method started at $x_0$, that is, any method that in round $t=0,1,\dots$ calls the LMO of $C$ with $\nabla f(z_t)$
 for a point $z_t$ of the convex hull $\mathcal K_t$ of $x_0$ and the answers of the previous
 rounds. Then $f(z)-\min_Cf\ge31/(256(t+1))$ for all $t\ge0$ and $z\in\mathcal K_t$.
\item If moreover $z_0=x_0$ and $f(z_{t+1})\le f(z_t+\frac2{t+2}(s_t-z_t))$ for the answer $s_t$
 of round $t$, as for Frank--Wolfe with step sizes $2/(t+2)$, then for all $t\ge1$
 \[
  \frac{31}{256(t+1)}\;\le\;f(z_t)-\min_Cf\;\le\;\frac{2(15/4)^2}{t+2}.
 \]
\end{enumerate}
\end{theorem}
\begin{proofsketch}
Write horizontal parts as $S\xi$. If the query point has $\xi$ supported on $[J]$ with mass
$\sum_i\xi_i\in[\frac12,\frac54]$, then the rows of the Gram matrix beyond $J$ are constant on
the support, so in the direction of its gradient index $J+1$ beats all later indices, and the
answer has $\xi$ supported on $[J+1]$ with mass in $[\frac12,1]$. By induction the hull after $t$
calls consists of such points with support in $[t+1]$, and Lemma~\ref{lem:hilbert-set}\,(ii) gives
$f(z)-\min_Cf\ge\frac12\norm{S\xi}^2\ge\frac{31}{64}\norm\xi^2\ge\frac{31}{256(t+1)}$. The upper
bound is the usual Frank--Wolfe bound with the diameter $15/4$.
\end{proofsketch}

\begin{remark}[Why finite dimension differs\leanverified]\label{rem:exposure}
The points $p_j=(d_j,\tau_j^2-1)$ lie in $C$, their height deficits $\tau_j^2$ tend to zero, and
$\norm{p_j-z^\star}\ge9/10$. So $z^\star$, the unique minimizer of the height and hence of
$\ip{\nabla f(z^\star)}\cdot$, is not strongly exposed: this is the step of
Proposition~\ref{prop:little-o} that fails, and in particular $C$ is not norm compact. A single arbitrary
query in the direction $(0,1)$ returns $z^\star$.
\end{remark}

\subsection{Arbitrary queries}

Hull-based methods query only gradients at points of the convex hull of their previous
answers. Methods with arbitrary queries may ask any direction, and then the set must depend on
the method; otherwise, as we have seen on the set of Theorem~\ref{thm:hilbert}, one query in the vertical
direction returns the minimizer (Remark~\ref{rem:exposure}). For deterministic methods one
set per method suffices for every horizon; the methods are anytime methods
(Definition~\ref{def:method}).

\begin{theorem}[Deterministic methods\leanverified]\label{thm:anytime}
Let $\mathcal H_0$ be an infinite-dimensional real Hilbert space, $v\in\mathcal H_0$ a unit
vector, $p=(0,-32)$ and $f(w)=\frac12\norm{w-p}^2$ on $\mathcal H_0\times\R$. For every
deterministic anytime method there are a weakly compact convex set $C\subseteq\mathcal
H_0\times\R$ of diameter at most $15/4$ that is $1$-smooth and has a unique minimizer for every
nonzero linear function, an LMO of $C$ that returns $x_0=(\sqrt2\,v,0)\in C$ for $q=0$, and a
minimizer $z^\star$ of $f$ on $C$, such that, run with this LMO, the output $w_T$ after $T$
calls satisfies, for every $T\ge0$,
\[
 w_T\notin C\qquad\text{or}\qquad f(w_T)-f(z^\star)\;\ge\;\frac1{4(T+2)} .
\]
\end{theorem}
\begin{proofsketch}
The set is $\bigcup_\lambda\bigl((h+\sum_i\lambda_ie_i,1)+(1+\sum_ib_i(\lambda_i-\lambda_i^2/8))B_1\bigr)$
over $\lambda\in\Delta$, with $b_i=2^{-i}/512$ ($i\ge0$) and an orthonormal system $h,e_0,e_1,\dots$
built along the run: the $t$-th new direction is chosen orthogonal to the horizontal parts of
all queries and outputs so far, and the system is arranged so that the start, known to the
method, is $(\sqrt2\,v,0)$. A query orthogonal to all $e_i$ with $i\ge k$ has its unique answer supported
on the first $k+1$ indices, since the decreasing radius increments $b_i$ favor the first unused
index and
the concave penalty makes the maximizer unique. Hence after $T$ calls the horizontal part of
the output is orthogonal to the witness $h-\sum_{i\le T}e_i$, of squared norm $T+2$, while the
horizontal part of the minimizer of $f$ has inner product at least $4/5$ with it. So a feasible
output has distance at least $(4/5)/\sqrt{T+2}$ from $z^\star$, and as $z^\star$ minimizes the
$1$-strongly convex $f$ on $C$, $f(w_T)-f(z^\star)\ge\frac12\norm{w_T-z^\star}^2\ge\frac{8}{25(T+2)}\ge\frac1{4(T+2)}$,
with one minimizer for all $T$.
\end{proofsketch}

For randomized methods we obtain two results, both up to logarithmic factors: a single
probability distribution over sets that is hard for every randomized method at every horizon;
and, for each randomized method, a single set that depends on the method but not on the outcome
$\omega$ of its random experiment, on which the method, with probability at least $1/5$ over $\omega$,
returns no $\eps$-solution at infinitely many horizons. As the set is now random, measurability does not suffice,
and we use jointly measurable methods (Definition~\ref{def:method}).\footnote{For $C_\theta=[\theta,\theta+2]$ with $\theta$ uniform on $[0,1]$, a single positive query
returns $\theta$. The deterministic method that then returns $5$ if $\theta$ lies in a fixed
non-measurable set and $1$ otherwise is measurable in the sense of Definition~\ref{def:method}, but not
jointly measurable, and its
infeasibility event is not measurable.} Let
\[
 g(T)=\frac1{(T+2)\log(T+3)\,\bigl(\log\log(T+e^e)\bigr)^2}.
\]

\begin{theorem}[Randomized methods\leanverified]\label{thm:random}
There are a probability distribution on weakly compact convex $1$-smooth sets
$C\subseteq\ell_2\times\R$ of diameter at most $15/4$ with unique linear minimizers, all in the
support of the distribution, and a point $x_0$ in all of them, such that, with $c=7\cdot10^{-5}$,
$f(w)=\frac12\norm{w-p}^2$, $p=(0,-32)$, and LMOs returning $x_0$ for $q=0$:
\begin{enumerate}
\item every jointly measurable randomized method requiring at most $T$ oracle calls returns no
 $\eps$-solution for $\eps=c\,g(T)$ on the random set $C$ with probability at least $1/5$, over $C$ and $\omega$;
\item for every jointly measurable randomized anytime method there is a set $C$, depending on
 the method but not on $\omega$, on which it returns no $\eps$-solution for $\eps=c\,g(T)$ for infinitely many $T$
 with probability at least $1/5$.
\end{enumerate}
\end{theorem}
\begin{proofsketch}
Block $j\ge0$ of $\ell_2$, of dimension $256k_j$ with $k_j=2^j$, carries a hard set of
Section~\ref{sec:vertex} scaled by $\zeta_j^{1/2}$, where $\zeta_j=1/\log(j+4)-1/\log(j+5)$ has a telescoping sum,
and all blocks are rounded once by a common unit ball. Instead of the vertices, the radius is
perturbed by independent random grid amounts per index; this keeps all answers unique, and each
answer reveals one index per block, also for directions that vanish on a block. For the least
$j$ with $k_j\ge T+1$, block $j$ alone forces infeasibility or primal gap $\zeta_j/(400k_j)\ge c\,g(T)$
with probability at least $\frac14-\frac1{100}$: by Lemma~\ref{lem:posterior} and Step~2 of the proof of
Theorem~\ref{thm:vertex} the scaled output block misses the barycenter of the hidden index set with
probability at least $\frac14$, and a coupling as in Section~\ref{sec:proof-main} loses at most
$\frac1{100}$. Reverse Fatou's lemma turns~(a) into~(b).
\end{proofsketch}

\section{Comparison and discussion}\label{sec:discussion}

For LMO-span methods on a $\beta$-smooth rounded simplex in dimension $d>T$,
\citet[Theorem~4.1]{grimmer2026lower} prove for a $1$-smooth, $1$-strongly convex quadratic the
lower bound
\[
 \max\left\{0,\ \sqrt{\frac{d-T}2}\left(\frac{\diam C-2/\beta}{\sqrt2\,d}+\frac1{\beta\sqrt d}\right)-\frac1{\sqrt2\,\beta}\right\}^2
\]
after $T$ steps. With $\delta=\diam C-2/\beta$, the expression inside the maximum equals the
quantity $A_d$ of the next proposition (our notation), so their bound is $\max\{0,A_d\}^2$.

\begin{proposition}[The bound of Grimmer and Liu at fixed parameters\leanverified]\label{prop:gl}
For $d>T>0$, $\delta\ge0$ and $\beta>0$,
\[
 A_d=\frac{\delta}{2\sqrt d}\sqrt{1-T/d}-\frac{1-\sqrt{1-T/d}}{\sqrt2\,\beta}
 \;\le\;\frac{\beta\delta^2}{4\sqrt2\,T}.
\]
\end{proposition}
\begin{proof}
Put $v=T/d\in(0,1)$ and $y=1/\sqrt d$, so that $v=Ty^2$. Since $\sqrt{1-v}\le1$ and
$1-\sqrt{1-v}\ge v/2$,
\[
 A_d\;\le\;\frac{\delta y}2-\frac{v}{2\sqrt2\,\beta}\;=\;\frac{\delta y}2-\frac{Ty^2}{2\sqrt2\,\beta}.
\]
The right-hand side is a concave quadratic in $y$; it is maximal at $y=\beta\delta/(\sqrt2\,T)$,
where it equals $\beta\delta^2/(4\sqrt2\,T)$.
\end{proof}

Their bound is thus at most $\beta^2\delta^4/(32T^2)$ in every dimension, of order $1/T^2$ at
fixed $\beta$ and $D$, while Theorem~\ref{thm:main} gives order $1/T$ for unrestricted randomized
methods and arbitrary outputs. Lower bounds on the simplex
\citep{jaggi2013revisiting,lan2013complexity} rest on the sparsity of iterates; here outputs are arbitrary. On strongly
convex sets the order is $1/T^2$ \citep{pokutta2026lower,halbey2026lower,HRP2026CurvatureFWLowerBounds}.

The sets of Theorem~\ref{thm:main} are smooth and strictly convex, yet the worst case is of order
$1/T$. Linear-oracle methods are accelerated by uniform convexity of the set
\citep{garber2015faster,KAP2020}, as in Theorem~\ref{thm:tangent}, or by a local version of it
\citep{P2026FrankWolfe}, with constants that depend on the modulus of convexity. The hard sets of
Section~\ref{sec:vertex} have flat faces, so they are not uniformly convex
(Proposition~\ref{prop:hard}), and with the vertex LMO they are not locally dually sharp
(Proposition~\ref{prop:optimum}). The sets of Theorem~\ref{thm:main} are strongly convex for each
$T$, but only with moduli of uniform convexity at most $\beta/(10^8(T+1))$
(Proposition~\ref{prop:strongly-convex}). So the faster rates hold on each of these sets, but with
constants that are not uniform in $T$.
Proposition~\ref{prop:little-o} and Theorem~\ref{thm:hilbert} separate the per-horizon from the single-set view: on a
single finite-dimensional set with unique minimizers and nonzero optimal gradient Frank--Wolfe
is $o(1/t)$, but not uniformly over the class, and in infinite dimension not at all;
Theorems~\ref{thm:anytime} and~\ref{thm:random} extend the single-set view to arbitrary queries. The oracle matters, too: inner-ball and
parallel-body descriptions of smooth sets are classical
\citep{schneider2014convex,liu2023gauges,luner2026performance}, gauge, normal or membership
information exploits smoothness \citep{liu2023gauges,hazan2026linearly}, but through an LMO it
is visible only as the known ball summand of Lemma~\ref{lem:rolling}.

Open questions remain: the constants $1/64$ and $2$ are far apart; the lower bounds need
dimension of order $T$, as usual \citep{nemirovski1983problem}, and the regime of fixed
dimension is open; and for randomized methods with arbitrary queries we do not know a single
set that is hard at every horizon at the exact rate $1/T$: Theorem~\ref{thm:random} loses logarithmic
factors and holds at infinitely many horizons.

\section*{AI use statement}
The construction was proposed by the author. The proofs and the original Lean formalization
were produced with the Agentic Researcher framework \citep{zimmer2026agentic}. The port to the
Lemmata Repository, the generalized and additional statements, the figures and their scripts, and
this manuscript were developed in a dialog with Claude (Anthropic) under the author's direction,
editing and review; the Hilbert-space results of Section~\ref{sec:hilbert} and their formalization were developed
with OpenAI Codex. The checks include machine-checked proofs for all results marked \leanverified,
exact-arithmetic verification of constants, and numerical assertions in the plotting scripts.
All claims have been verified by the author, who takes responsibility for the final content,
including text, claims and artifacts produced with the aid of generative AI.

\ificmlbuild\else
\section*{Acknowledgements}
This research was conducted in the context of the Agentic AI in Mathematics (EF-LiOpt-3)
project at the Berlin Mathematics Research Center MATH$^+$ (EXC-2046/2, project ID 390685689),
funded by the Deutsche Forschungsgemeinschaft (DFG, German Research Foundation) under
Germany's Excellence Strategy. This work has been partially supported by OpenAI's Academic
Researchers program.
\fi

\ificmlbuild
\section*{Impact Statement}
This paper presents work whose goal is to advance the field of Machine Learning. There are many
potential societal consequences of our work, none which we feel must be specifically highlighted
here.

\bibliographystyle{icml2027}
\else
\bibliographystyle{plainnat}
\fi
{\IfNotIcmlBuild{\small\setlength{\bibsep}{1pt plus 1pt}}
 }

\clearpage
\appendix
\IfIcmlBuild{\onecolumn}
\section{Proof of \texorpdfstring{Theorem~\ref{thm:main}}{Theorem 3.2}}\label{app:main}

We use the reduction to $r=L=1$, $f(x,y)=\frac12(\norm x^2+y^2)+32y$, and the notation $k,N,n,A_0,A$
of Section~\ref{sec:proof-main}; note $x_0=r(Ae_1,0)$.

\begin{proof}[Proof of Theorem~\ref{thm:main}\,(a)]
The proof follows the sketch in Section~\ref{sec:proof-main}. We perturb the vertices by a random
table~$Z$ (\emph{The sets}); show that each answer then reveals one index of $I$ unless a comparison
is ambiguous (\emph{Each answer reveals one index}); simulate the given method by a method that
learns indices of $I$ (\emph{Simulation}); show that ambiguous comparisons are rare (\emph{Ambiguous
answers are rare}); compare the perturbed with the unperturbed sets (\emph{Transfer}); and average
(\emph{Averaging}).

\emph{The sets.} Let $\eta=A^2/(10^6k(A+36)^2)$, $m_0=0$, $m_{t+1}=1+nm_t$ (so $m_T$ bounds the
number of nodes of a tree of depth $T$ and branching $n$), $Q=1+100m_TN^2n$, $\kappa=\eta/(nQ)$ and
$\sigma=\kappa/(8n)$: $\eta$ bounds the perturbation, $Q$ is the number of grid values per coordinate,
$\kappa$ the grid step, and $\sigma$ the thickening. A \emph{table} $Z$ assigns to every index $j\ne1$ a vector $Z_j\in\R^n$ with
independent coordinates, uniform on $\{0,\kappa,\dots,(Q-1)\kappa\}$, and $Z_1=0$; then
$\norm{Z_j}\le\eta$. For $I\ni1$ let $v_j=(Ae_j,1)+Z_j$ and
\[
 C_{I,Z}=\bigcup_{\lambda\in\Delta_I}\Bigl(\sum_{j\in I}\lambda_jv_j+\bigl(1+\sigma(1-\norm\lambda^2)\bigr)B_1\Bigr).
\]
As $1-\norm\lambda^2$ is concave, $C_{I,Z}=K_{I,Z}+B_1$ with $K_{I,Z}$ convex, so $C_{I,Z}$ is a
$1$-smooth convex body (as in the proof of Proposition~\ref{prop:hard}) containing $v_1+B_1\ni(Ae_1,0)$.
For a unit $q$, minimizing $\ip q\cdot$ over $C_{I,Z}$ amounts to minimizing the strictly
convex function $\sum_{j\in I}\lambda_j\ip q{v_j}+\sigma\norm\lambda^2$ over $\Delta_I$, so the
minimizer is unique. With the same weights, centers move by at most $\eta$ and radii by at
most $\sigma$ relative to $C^0_I=\conv\{(Ae_j,1):j\in I\}+B_1$, so the Hausdorff distance is at
most $\vartheta:=2\eta$, and $\diam C_{I,Z}\le2+\sqrt2A+2\vartheta\le\beta D$.

\emph{Each answer reveals one index.} For a unit $q$ let $j^\star$ be the first index of $I$ in
the order of $(\ip q{v_j},j)$. If $j\in I$ has $\ip q{v_j}>\ip q{v_{j^\star}}+2\sigma$, the
minimizing weights vanish at $j$: moving weight from $j$ to $j^\star$ changes the objective at
rate $\ip q{v_{j^\star}}+2\sigma\lambda_{j^\star}-\ip q{v_j}-2\sigma\lambda_j<0$.

\emph{Simulation.} Fix $\omega$ and $Z$. The method $M_Z$ learns with each call the first index
of $I$ in an order of its choice and runs the given method, keeping the revealed indices $R$,
initially $\{1\}$. When the given method asks $q\ne0$, $M_Z$ asks for the first index $j^\star$
of $I$ in the order of $(\ip{q/\norm q}{v_j},j)$, adds it to $R$, and passes on the minimizer
of $\ip q\cdot$ over $C_{R,Z}$; for $q=0$ it passes on $x_0$. When the given method returns $w$,
$M_Z$ returns $\widehat w=(w+\vartheta v_1)/(1+\vartheta)$. The orders depend on $Z$, which is
independent of $I$, and the passed answers on the revealed vertices only. If at every nonzero
query all $j\in I\setminus R$ satisfy $\ip q{v_j}>\ip q{v_{j^\star}}+2\sigma$, the minimizers over
$C_{I,Z}$ and $C_{R,Z}$ coincide, and the given method run with $\LMO_{C_{I,Z},x_0}$ produces
the simulated transcript and output.

\emph{Ambiguous answers are rare.} Call $Z$ \emph{good} if at every node of the decision tree
of $M_Z$, with unit query $q$, all pairs $i\ne j$ with $j$ not yet revealed satisfy
$|\ip q{v_i}-\ip q{v_j}|>2\sigma$; then the condition above holds for every $I$. Fix a node, a
pair $i\ne j$ with $j\ne1$ unrevealed, and a coordinate $c$ with $|q_c|\ge1/n$. The query
depends only on $\omega$ and the revealed vertices, so given all entries of $Z$ but $Z_{j,c}$,
$\ip q{v_i}$ is fixed, while one grid step of $Z_{j,c}$ changes $\ip q{v_j}$ by at least
$\kappa/n=8\sigma$. So at most one of the $Q$ values of $Z_{j,c}$ is ambiguous, and a union
bound over at most $m_T$ nodes, $N^2$ pairs and $n$ coordinates shows that $Z$ is bad with
probability at most $m_TN^2n/Q<1/100$.

\emph{Transfer.} Let $Z$ be good and $w\in C_{I,Z}$. For $v\in C^0_I$ with $\norm{w-v}\le\vartheta$,
$\widehat w$ is a convex combination of $v$ and $v_1+(w-v)/\vartheta\in C^0_I$, so
$\widehat w\in C^0_I$, and $\norm{w-\widehat w}\le\vartheta(2A+4)$. Both sets lie in the ball of
radius $A+3$, where $f$ is $(A+35)$-Lipschitz, and their minima of $f$ differ by at most
$(A+35)\vartheta$. Hence $f(w)-\min_{C_{I,Z}}f\ge f(\widehat w)-\min_{C^0_I}f-(A+35)(2A+5)\vartheta$,
and $(A+35)(2A+5)\vartheta\le4A^2/(10^6k)$.

\emph{Averaging.} For every $\omega$ and good $Z$, Steps~1--2 of the proof of
Theorem~\ref{thm:vertex}, applied to $M_Z$ on the sets $C^0_I$ (scale $A$, radius $1$), show that
$\widehat w$ is no $\eps$-solution on $C^0_I$ for $\eps=A^2/(16k)$ for at least a quarter of the index sets;
then $w$ is no $\eps$-solution on $C_{I,Z}$ for $\eps=A^2/(16k)-4A^2/(10^6k)\ge A_0^2/(32k)$. So for every
$\omega$, on a fraction at least $\frac14\cdot\frac{99}{100}\ge\frac15$ of the pairs $(I,Z)$, the
method returns no $\eps$-solution, and averaging over $\omega$ selects one pair. The set $rC_{I,Z}$ lies in
$\cls(\beta,D,x_0)$, and after scaling $\eps$ becomes $Lr^2A_0^2/(32k)=L(D-2/\beta)^2/(64(T+1))$.
\end{proof}

\begin{proof}[Proof of Theorem~\ref{thm:main}\,(b)]
The sets $rC_{I,Z}$ lie in the class, so the given method returns points of them, and
$w\mapsto\widehat w$ moves these into the unperturbed sets. For good $Z$, $M_Z$ thus has feasible
outputs, and the argument of Theorem~\ref{thm:variants}\,(a) (Appendix~\ref{app:variants}) gives no
$\eps$-solution on a fraction at least $97/128$ of the index sets; $\frac{97}{128}\cdot\frac{99}{100}\ge\frac34$. The primal
gap is nonnegative and bounded, which gives the expectation.
\end{proof}

\section{Proof of \texorpdfstring{Theorem~\ref{thm:variants}}{Theorem 4.6}}\label{app:variants}

We use the notation of the proof of Theorem~\ref{thm:vertex}: $k=T+1$, $d=8k$, the weights
$\theta=\theta(w)$, and, for fixed $\omega$ and uniform $I$, the sets $P,U$, the number $h$ and
the membership probabilities $p_i$ given the final answers.

\begin{proof}[Proof of (a)]
Fix $\omega$. \emph{One light instance suffices.} All instances consistent with the final
answers give the same output $w$, which lies in each of their sets. If one of them, $J$, has
$\sum_{i\in J}\theta_i\le3/4$, Step~1 applied to $J$ gives $f(w)\ge La^2/(8k)$, hence
$f(w)-\min_{C_I}f\ge La^2/(16k)$ for every consistent $I$. \emph{Many undecided indices leave a
light instance.} If $|U|\ge127k$, then $\sum_ip_i^2=|P|+h^2/|U|\le k+64k/127<9k/4$, and
$p_i\theta_i\le p_i^2/(6k)+\frac{3k}2\theta_i^2$ shows that the consistent instances have
average weight at most $3/4$. \emph{Few indices are excluded.} A decided index lies in $I$, at
most $k$ times, or is excluded; so at most $128k$ exclusions leave $|U|\ge127k$. Summed over the
$\binom{|U|}h$ consistent instances, a run with $t$ calls left excludes at most
$t\binom{|U|}{h+1}$ indices: at each decided index both binomial coefficients split by
Pascal's rule, an exclusion costs no call, and an index found in $I$ ends the call. At the
start this is $31T$ per instance, so by Markov's inequality more than $128k$ exclusions have
probability at most $31T/(128k)<31/128$, and the output is no $\eps$-solution with probability at least
$97/128\ge3/4$. Averaging over $\omega$ selects one $I$; the primal gap of a feasible point is
nonnegative and bounded, which gives the expectation.
\end{proof}

\begin{proof}[Proof of (b)]
Answer every query $q\ne0$ as if the first index in its order belonged to $I$, and record it.
With $P$ the recorded indices and $1$, $|P|\le k$, and every index set containing $P$ produces
exactly these answers and the same output. Complete $P$ by indices chosen uniformly among the
remaining $N-|P|\ge63k$. Then $\sum_ip_i^2\le k+64k/63<9k/4$, so as in~(a) some completion $I$ has
$\sum_{i\in I}\theta_i\le3/4$, and by Step~1 the output is no $\eps$-solution on it.
\end{proof}

\begin{proof}[Proof of (c)]
Let $v$ be the point of $C_I$ nearest to $w$ and $s=\norm{w-v}$. On $C_I$ the gradient of $f$ has
norm at most $L(a+34r)$, so $f(w)\ge f(v)-L(a+34r)s$, and $\rho\ge L(2\sqrt2a+36r)\ge L(2a+34r)$
gives $\err_I(w)\ge f(v)-\min_{C_I}f+Las$. By scaling, $a=L=1$; the bound $\eps/8$ is then
$1/(128k)$. On the event $\sum_{i\in I}\theta_i\le3/4$, of probability at least $1/4$ by Step~2:
if $s\ge1/(128k)$ we are done. Otherwise, if $\norm x^2>1/(4k)$ for the horizontal part $x$ of
$w$, the horizontal part of $v$ has norm at least $1/(2\sqrt k)-1/(128k)\ge3/(8\sqrt k)$, and
$f(v)\ge9/(128k)$; if not, the horizontal coordinates of $v$ on $I$ sum to at most
$3/4+\sqrt{8k}\,s\le4/5$, and the argument of Step~1 gives $f(v)\ge2/(25k)$. Subtracting
$\min_{C_I}f\le1/(16k)$ leaves at least $1/(128k)$. Averaging over $\omega$ selects one $I$, and
the mean is at least a quarter of this bound.
\end{proof}

\PrintDeferredProofs

\end{document}